\documentclass[
a4paper]{amsart}

\usepackage{dropsty}
\usepackage{fixme}

\usepackage{setspace}

\title{Group action combinatorics}
\author{Brendan Murphy}
\date{\today}

\begin{document}

\begin{abstract}
  This paper generalizes the basic notions of additive and multiplicative combinatorics to the setting of \emph{group actions}:
if $G$ is a group acting on a set $X$, and we have subsets $A\subseteq G$ and $Y\subseteq X$ such that the set of pairs $g\cdot y$ with $g\in A,y\in Y$ is not much larger than $Y$, what structure must $A$ and $Y$ have?
  Briefly, what is the structure of sets with small \emph{image set}?
  
  In this setting, we develop analogs of Ruzsa's triangle inequality, covering theorems, multiplicative energy, and the \bsg{} theorem.
Approximate stabilizers, or \emph{symmetry sets}, play an important role.
  
  While our focus is on presenting a general theory, we answer the inverse image set question in some special cases.
  To do so, we combine the group action version of the \bsg{} theorem with structure theorems for approximate groups and bounds for the sizes of symmetry sets.
\end{abstract}

\maketitle 
\tableofcontents

\addtocontents{toc}{\protect\setcounter{tocdepth}{1}}
\section{Introduction}

In the series of papers \cite{elekes1997linear,elekes1998linear,elekes1999linear}, Elekes studied a \emph{non-commutative}\/ version of Freiman's theorem.
Namely, if $L$ is a set of $N$ affine transformations $\ell(x)= mx+b$, for some collection of pairs of real numbers $(m,b)$, and $A$ is a set of $N$ real numbers, Elekes studied the \emph{image set}\/ $L(A)=\{\ell(a)\colon \ell\in L, a\in A\}$ and asked: ``What structure must $L$ and $A$ have if $|L(A)|\leq KN$?''

If $L$ is a set of translations (that is, $m=1$ for all $\ell(x)=mx+b$ in $L$), then Freiman's theorem implies that $A$ and the set of $b$'s corresponding to lines in $L$ are both contained in \emph{generalized arithmetic progressions}.
Similarly, if $L$ is a set of dilations ($b=0$ for all $\ell(x)=mx+b$ in $L$), the $A$ and the set of $m$'s corresponding to the lines in $L$ are both contained in \emph{generalized geometric progressions}.
Elekes proved that any set of affine transformations $L$ must contain a large subset with one of these two structures.
In later work, Elekes and Kiraly generalized this result to linear fractional transformations \cite{elekes2001combinatorics}, and conjectured analogous theorems for actions of higher dimensional matrix groups.

This paper expands on Elekes' framework by establishing group action analogs of tools from multiplicative combinatorics.
Further, we prove a group action version of the \emph{asymmetric \bsg{} theorem}, which allows us to extend the results of \cite{elekes1997linear,elekes1998linear,elekes1999linear,elekes2001combinatorics,elekes2002versus} to any action of an algebraic group.
In fact, we can prove such theorems with much weaker hypotheses; the
paper \cite{murphy2017upper} uses this method to improve work of Croot and others \cite{borenstein2010lines,amirkhanyan2017lines}, which extended Elekes' work on rich lines in grids.

For motivation, we consider a general conjecture on approximate group actions in the next section.

\subsection*{Overview of the problem and methods involved}

Let $G$ be a group acting on a set $X$, let $A$ be a subset of $G$, and let $Y$ be a subset of $X$.
We use $A(Y)$ to denote the set of points $a(y)\in X$ with $a\in A$ and $y\in Y$; this is the \emph{image set}\/ of $Y$ under $A$.
The following theorem characterizes when $A(Y)$ is not larger than $Y$.
\begin{thm}
  \label{thm:1}
Let $G$ be a group acting on the set $X$, let $A$ be a subset of $G$ and let $Y$ be a finite subset of $X$.
Suppose that $|A(Y)|=|Y|$.
Then $A^{-1}A$ generates a subgroup $H$ of the stabilizer $\stab(Y)$ and $Y$ is a union of $H$-orbits.
\end{thm}
Here, $\stab(Y)$ denotes the stabilizer of the set $Y$ under the induced action of $G$ on subsets of $X$.
\begin{proof}
  Since $e\in a^{-1}A$, we have $Y\subseteq a^{-1}A(Y)$.
But $|Y|=|A(Y)|=|(a^{-1}A)(Y)|$, so $a^{-1}A(Y)=Y$.
Since $a\in A$ was arbitrary, it follows that $A^{-1}A(Y)=Y$, hence $A^{-1}A\subseteq \stab(Y)$.
Let $H$ be the subgroup of $\stab(Y)$ generated by $A^{-1}A$.
Since $H$ acts on $Y$, it follows that $Y$ is a union of $H$-orbits, as claimed.
\end{proof}

The starting point for this paper is relaxing the conditions in Theorem~\ref{thm:1}.
Suppose instead that $|A(Y)|\leq K|Y|$ for some parameter $K\geq 1$.
Is it still true that $Y$ is \emph{approximately}\/ a union of orbits of a subgroup of $G$ generated by $A$ or $A^{-1}A$?
\begin{conj}
  \label{conj:main}
Let $G$ be a group acting on a set $X$, let $A$ be a subset of $G$ and let $Y$ be a finite subset of $X$.
Suppose that $|A(Y)|\leq K|Y|$ for some $K\geq 1$.
Then there is a constant $C>0$, a subset $B\subseteq G$, a subgroup $H\leq G$, and a finite subset $Z\subseteq X$ such that $|B|\ll K^C$, $A\subseteq BH$, and $|H(Z)\cap Y|\gg K^{-C}|Y|$.
\end{conj}

The proof of Theorem~\ref{thm:1} suggests introducing an approximate analog of the stabilizer of $Y$.
For $0<\alpha\leq 1$, we let $\sym_\alpha(Y)$ be the set of $g\in G$ such that $|Y\cap gY|\geq\alpha|Y|$; this is a \emph{symmetry set}\/ of $Y$.
As before, for any $a\in A$, we have $Y\subseteq a^{-1}A(Y)$.
Setting $Y'=a^{-1}A(Y)$, it follows that for any $g\in a^{-1}A$, we have $gY\subseteq Y'$, hence $|Y'\cap gY'|\geq |gY|\geq \frac 1K |Y'|$, so $a^{-1}A\subseteq\sym_{\alpha}(Y')$ with $\alpha=\frac 1K$.
This reduces the problem to studying symmetry sets.

We want to show that symmetry sets behave like groups.
As a first step, we show that symmetry sets have \emph{weak multiplicative closure}\/ (Proposition~\ref{prop:3}).
Using an iteration scheme we prove that $\sym_\alpha(Y)$ is controlled by an \emph{approximate group}.
To close the iteration, we need bounds for $|\sym_\alpha(Y)|$.
To find more precise structure, we need \emph{structure theorems}\/ (or \emph{product theorem}) for approximate groups.

These two ingredients, symmetry set bounds and product theorems, limit what we can prove.
Still, our method is more flexible than Elekes', and we can prove Conjecture~\ref{conj:main} in some cases.

\subsection*{Related work}

This work is inspired by Elekes' work on rich affine and linear fractional transformations, as well as his approach to the Erd\H{o}s distance problem \cite{elekes1997linear,elekes1998linear,elekes1999linear,elekes2002versus,elekes2001combinatorics,elekes2011incidences,elekes2011dimension}.
My work on rich lines in grids \cite{murphy2017upper}, following that of Elekes and subsequent work \cite{borenstein2010lines,amirkhanyan2017lines}, can be read in parallel with this work.
Michael McGee informed me that the approach of this paper is similar to Bourgain's proof of an incidence theorem for modular hyperbolas \cite{bourgain2012modular}; comparing Bourgain's proof to the proof I give in Section~\ref{sec:low-dimens-exampl} may be instructive.
Finally, Harald Helfgott has emphasised that the force behind the sum-product problem results on growth in groups is the \emph{tensions between two group actions}; this philosophy was another inspiration for this paper.

\subsection*{Organization}

The rest of the paper is organized as follows:
\begin{itemize}
\item Section~\ref{sec:group-acti-comb} discusses the basics of group action combinatorics.
\item Section~\ref{sec:group-action-bsg} contains the statement and proof of the group action version of the \bsg{} theorem; this generalizes the asymmetric \bsg{} theorem from additive combinatorics.
The end of the section contains an application of the group action \bsg{} theorem to ``nearly free'' actions; in particular, this generalizes the asymmetric \bsg{} to non-commutative groups.
\item Section~\ref{sec:applications} contains some applications of the general theory.
The first application is an alternate proof of Bourgain's incidence theorem for hyperbolas \cite{bourgain2012modular}
The second application is a generalization of results of Elekes and Kiraly to any matrix group over a field of characteristic zero.
Roughly, this result is that if $A\subseteq GL_n(k)$, where $\mathrm{char}(k)=0$, $Y\subseteq k^n$ is finite, and $|Y\cap gY|\geq \alpha |Y|$ for all $g\in A$, then either $Y$ is (mostly) contained in a hyperplane or $A$ is (mostly) contained in a coset of a nilpotent subgroup of $GL_n(k)$.
\end{itemize}

\subsection*{Notation}

We use standard asymptotic notation: for positive functions $f$ and $g$, we write $f\ll g$ if there is a constant $C>0$ so that $f\leq Cg$; similarly $f\gg g$ means $g\ll f$; we also write $f=O(g)$ if $f\ll g$ and $f=\Omega(g)$ if $f\gg g$.
If we use a subscript, say $f\ll_r g$, then the implicit constant depends on $r$: $f\leq C(r)g$; thus $O_r(1)$ means a positive (unspecified) constant depending on $r$.

If $G$ is a group acting on a set $X$, we write $G\actson X$.
If $G\actson X$ and $Y\subseteq X$, we use $\stab(Y)$ to denote the \emph{set-wise stabilizer}\/ of $Y$.
That is, $\stab(Y)=\{g\in G\colon gY=Y\}$.
The \emph{point-wise stabilizer}\/ of $Y$ is $\bigcap_{y\in Y}\stab(y)$.
For $g\in G$, we use $\fix(g)$ to denote the set of \emph{fixed points}\/ of $g$:
\[
\fix(g)=\{x\in X\colon gx=x\}.
\]

For a subset $A$ of a group $G$, we define the product set $AA=\{aa'\colon a,a'\in A\}$ and use $A^k$ to denote the $k$-fold product of $A$ with itself.
We use $A^{-1}$ to denote the set of inverses of elements of $A$.
If $A=A^{-1}$, we say that $A$ is \emph{symmetric}.
It is often useful to assume that a set is symmetric, so we use the notation $A_{(k)}=(A\cup A^{-1}\cup\{e\})^k$ to denote $k$-fold products of the symmetrization of $A$.
We also use exponents to denote Cartesian products, but typically for sets $Y\subseteq X$, though of as sets of points in the space that $G$ is acting upon.

\subsection*{Acknowledgements}

The work of G. Elekes inspired this paper.
I would like to thank Giorgis Petridis, Ilya Shkredov, Misha Rudnev, Alex Iosevich,  Jonathan Pakianathan, and Harald Helfgott for helpful discussions and comments on drafts of this paper.
I would also like to acknowledge the Heilbronn Institute for Mathematical Research and the Leverhulme Trust for their support.

\section{Background}

\subsection{Group actions}

In this section, we review some basic ideas about groups actions.
A group $G$ is said to act on a set $X$ (on the left) if there is a map from $G\times X$ to $X$, which we will denote by juxtaposition $(g,x)\mapsto gx$, that satisfies
\[
ex=x\qquad\mbox{for all $x\in X$}\andd g(hx)=(gh)x\qquad\mbox{for all $g,h\in G,x\in X$}.
\]
We may refer to the elements of $g$ as ``transformations'' acting on the ``space'' $X$.
For an element $x\in X$, the \emph{stabilizer $\stab(x)$ of $x$}\/ is the set of transformations that fix $x$, and the \emph{orbit $Gx$ of $x$}\/ is the set of all points in $X$ that can be reached from $x$: $Gx:=\{gx\colon g\in G\}$.
The stabilizer of any point is a \emph{subgroup}\/ of the group $G$.
If $x'$ is in the orbit of $x$, then the stabilizer of $x'$ is \emph{conjugate}\/ to the stabilizer of $x$.

We say that an action is \emph{transitive}\/ if for any $x$ in $X$, the orbit of $x$ is the entire space: $Gx=X$.
Transitive actions can be described entirely within $G$.
For any $x$ in $X$, we define the \emph{orbit map}\/ $\phi_x\colon G\to X$ by $\phi_x(g)=gx$.
The orbit map is \emph{equivariant}\/, meaning that $\phi_x(gh)=g\phi_x(h)$ for any $g,h\in G$.
If $G$ acts transitively on $X$, then the orbit map is surjective and induces an equivariant bijection from $G/\stab(x)$ to $X$.
That is, the action of $G$ on $X$ and the action of $G$ on $G/\stab(x)$ are \emph{isomorphic}.
Letting $H=\stab(x)$, we have
\[
\tilde\phi_x\colon G/H\to X
\]
defined by $gH\mapsto gx$.
Note that the orbit map is not a group homomorphism unless $X$ is a group.

The above isomorphic leads to the \emph{orbit-stabilizer theorem}, which states that if a finite group $G$ acts transitively on a space $X$, then for any point $x\in X$ we have
\[
|G| = |Gx||\stab(x)|.
\]
Approximate version of this statement play an important role in approximate group theory.

Further terminology used in the paper includes:
\begin{itemize}
\item $n$-fold transitive, meaning that $G$ acts transitively on the set of $n$-tuples of distinct points of $X$
\item faithful, meaning that no element of $G$ besides the identity fixes every element of $X$
\item free, meaning that no element of $X$ has a non-trivial stabilizer (that is, no element is fixed by any transformation besides the identity)
\end{itemize}

Another concept that is important to this paper is the \emph{transporter}\/ from a point $x$ to a point $x'$:
\[
\trans(x,x'):=\{g\in G\colon gx=x'\}.
\]
The transporter from $x$ to $x'$ is non-empty only if $x'$ is in the orbit of $x$.
If $x'=g_0x$, then
\[
\trans(x,x')=g_0\stab(x)=\stab(x')g_0,
\]
so transporters are cosets of stabilizers.

Further references for group actions are \cite{artin1991algebra} and \cite{lang2002algebra}.
The concepts discussed here tend to behave well in situations with additional structure, for instance for algebraic or smooth actions \cite{humphreys1975linear,borel1991linear,pontryagin1966topological}.

\subsection{Basic combinatorics}

\begin{lem}[Popularity principle]
  \label{lem:Pop}
Given a set $X$ and a positive function of finite support $f\colon X\to \R_{\geq 0}$, let
\[
m=\frac 1{|\supp f|}\sum_x f(x).
\]
Fix $0 < \lambda < 1$ and let $P_\lambda = \{x\in X\colon f(x)\geq \lambda m\}$.
Then
\[
\sum_{x\in P_\lambda}f(x) \geq (1-\lambda)\sum_x f(x).
\]
\end{lem}

\begin{lem}[Cauchy-Schwarz intersection lemma]
\label{lem:CSint}
 Let \(S\) be a finite index set and let \(T_s\) be a family of subsets of
a set \(T\).
Then
\begin{equation*}%
  \left(\sum_{s\in S}|T_s| \right)^2\leq |T|\sum_{s,s'\in S}|T_s\cap T_{s'}|.
\end{equation*}%
Further, if there exists \(\delta > 0\) such that
\[
 \sum_{s\in S}|T_s| \geq \delta |S||T|,
\]
then there exists a subset \(P\subseteq S\times S\) such that
\begin{enumerate}
\item \( |T_s\cap T_{s'}|\geq \delta^2 |T| / 2 \) for all pairs \((s,s')\)
   in \(P\).
\item \( |P|\geq \delta^2 |S|^2 / 2 \).
\end{enumerate}
\end{lem}

\addtocontents{toc}{\protect\setcounter{tocdepth}{2}}
\section{Group action combinatorics}
\label{sec:group-acti-comb}

In this section, we generalize the basic theorems of additive and arithmetic combinatorics to the setting of group actions.
Throughout, $G$ denotes a group acting on a set $X$.
Two basic examples to keep in mind are the following.
\begin{ex}[Additive/multiplicative combinatorics]
  \label{ex:mult-comb}
In \emph{multiplicative combinatorics}, a group $G$ acts on itself by left translation; thus, $X=G$.
When $G$ is abelian and written additively, we say \emph{additive combinatorics}\/ instead.
\end{ex}
\begin{ex}[1-dimensional affine transformations]
  \label{ex:1D-aff}
Let $\F$ be a field, let $X=\F$, and let $G=\aff 1\F$ be the group of affine transformations of $X$.
An element of $G$ has the form $x\mapsto ax+b$, where $a,b\in\F$, $a\not=0$.
\end{ex}

This section is organized as follows:
\begin{itemize}
\item In Section~\ref{sec:image-sets}, we define \emph{image sets}; for $G\actson G$ by left translation, an image set is a product set.
\item In Section~\ref{sec:symmetry-sets}, we define \emph{symmetry
    sets}.
For an abelian group acting on itself by translation, symmetry sets are sets of popular differences.
For the example of $\aff 1\F\actson \F$, symmetry sets correspond to sets of lines that contain many points of a Cartesian product point set (``rich lines'').
\item In Section~\ref{sec:energy}, we define \emph{action energy}, which specializes to multiplicative energy when $G\actson G$ by left translation.
The paper \cite{yazici2015growth} applies bounds on the action energy of $\aff 1\F\actson\F$ to prove sum-product theorems.
\item In Section~\ref{sec:conversions}, we show how to convert between image sets, symmetry sets, and action energy; for $G\actson G$ by left translation, this corresponds to the well-known results that small product set implies large energy, large energy implies many popular ratios and small partial product sets, and conversely.
\end{itemize}

\subsection{Image sets}
\label{sec:image-sets}

If $A$ is a subset of $G$ and $Y$ is a subset of $X$, we define the \emph{image set} of \(Y\) under \(A\) by
\[
A(Y) = \{g(x)\colon g\in A, x\in Y\}.
\]
We also define a \emph{partial}\/ or \emph{statistical}\/ version of image sets.
For any subset \(E\)  of \(A\times Y\) we define the \emph{partial image set} of \(Y\) under \(A\) by
\[
A_E(Y) = \{g(x)\colon (g,x)\in E\}.
\]

Image sets unify several ideas from additive and arithmetic combinatorics.
\begin{repex}{ex:mult-comb}[Image sets in multiplicative combinatorics]
  For $G\actson G$ by left translation, the image set $A(Y)$ is the \emph{product set}\/ of $A$ and $Y$.
\end{repex}
\begin{repex}{ex:1D-aff}[Image sets for affine transformations]
  For affine transformations, if $Y$ is a subset of $X=\F$ and $L$ is
  a subset of $\aff 1\F$, then the image $L(Y)$ of $Y$ under $L$ was
  studied by several authors \cite{elekes1997linear,elekes1998linear,elekes1999linear,elekes2002versus,yazici2015growth} and has connections to the \emph{sum-product problem}.

For instance, if $L$ is the set of transformations of the form $x\mapsto a(x+b)$ where $a,b\in A$, then $L(A)=A(A+A)$, and if $L$ is the set of transformations of the form $x\mapsto a+bx$, then $L(A)=A+AA$.
\end{repex}

When $Y$ is a singleton, we say that $A(\{y\})$ is an \emph{approximate orbit}, and write $A(y)$ when no confusion may result.

An important example is the action of a group on the \emph{coset space}\/ of a subgroup.
\begin{ex}[Left multiplication on a coset space]
  \label{ex:coset-action}
Let $G$ be a group, $H$ a subgroup of $G$, and let $G$ act on the coset space $X=G/H$ by left multiplication.
Let $\pi\colon G\to G/H$ denote the canonical map ($\pi$ is a homomorphism only if $H$ is normal).

If $A\subseteq G$, then the approximate orbit $A(\{H\})$ is the image $\pi(A)$ of $A$ under the canonical map.
(We write $A(\{H\})$ to distinguish the approximate orbit from the product set $AH\subseteq G$, which is the image set $A(H)$ for $G\actson G$ by left-multiplication.)
\end{ex}

Now we generalize two basic results of multiplicative combinatorics to groups actions: Ruzsa's triangle inequality and Ruzsa's covering lemma.

\subsubsection{Triangle inequality for image sets}
Recall that Ruzsa's triangle inequality states that if $A,B,$ and $C$ are non-empty finite subsets of a group, then
\begin{equation}
  \label{eq:ruzsa-triangle}
  |AC^{-1}|\leq\frac{|AB^{-1}||BC^{-1}|}{|B|}.
\end{equation}
The group action version of \eqref{eq:ruzsa-triangle} is:
\begin{prop}[Ruzsa's triangle inequality]
\label{prop:1}
 Let \(A_1\) and \(A_2\) be non-empty finite subsets of \(G\).
Then for any finite subset $Y$ of $X$,
\[
|A_1||A_2(Y)| \leq |A_2A_1^{-1}||A_1(Y)|.
\]
\end{prop}
For multiplicative combinatorics ($G\actson G$), Proposition~\ref{prop:1} with $A_1=B, A_2=A,$ and $Y=C^{-1}$ is (\ref{eq:ruzsa-triangle}).
The proof of Proposition~\ref{prop:1} is essentially the same as the proof of (\ref{eq:ruzsa-triangle}).
\begin{proof}
To show that
\[
| A_1 ||A_2(Y)| \leq |A_2A_1^{-1}||A_1(Y)|,
\]
it suffices to find an injection
\[
\phi\colon A_1\times A_2(Y) \to A_2A_1^{-1}\times A_1(Y).
\]

For each element $x$ of $A_2(Y)$, we select a pair $(a_x,y_x)$ in $A_2\times Y$ such that $a_xy_x=x$.
Now for a pair $(a,x)$ in $A_1\times A_2(A)$, we define $\phi(a,x)$ by
 \[
\phi\colon (a,x)\mapsto (a_xa^{-1},ay_x).
\]
The image of $(a,x)$ is contained in $A_2A_1^{-1}\times A_1(A)$.
Further $\phi$ is injective, since we may recover the pre-image of any element $(a',x')$ in the image of $\phi$, by first finding $x=a'x'$ and then using $a_x$ fo solve for $a=(a')^{-1}a_x$.
\end{proof}

Proposition~\ref{prop:1} is familiar in additive and multiplicative combinatorics.
The following example, taken from the theory of approximate groups \cite{helfgott2015growth,breuillard2014brief}, involves the action of $G$ on the coset space $X:=G/H$.
\begin{cor}[Growth in a subgroup implies growth]
\label{cor:growth-in-a-subgroup}
Let $A$ and $B$ be non-empty finite subsets of a group $G$, let $H$ be a subgroup of $H$, and let $\pi\colon G\to G/H$ be the quotient map.
Then if $B\cap H$ is non-empty, we have
\[
|\pi(A)||B\cap H|\leq |AB|.
\]
\end{cor}
In particular, if $B=A^{N}$, then 
\[
|\pi(A)||A^N\cap H|\leq |A^{N+1}|,
\]
hence if $A\cap H$ grows, so does $A$.
\begin{proof}
Let $A$ and $B$ be subsets of $G$.

Applying Proposition~\ref{prop:1} with $A_1=B^{-1}\cap\stab(H)=B^{-1}\cap H$, $A_2=A$, and $Y=\{H\}$ yields
\[
|A(\{H\})||B^{-1}\cap H| \leq |A(B\cap H)||(B^{-1}\cap H)(\{H\})|\leq |AB|.
\]
As mentioned in Example~\ref{ex:coset-action}, $A(\{H\})$ is the image of $A$ under the canonical map $\pi\colon G\to G/H$.
Further, $|B^{-1}\cap H|=|B\cap H|$, so we have
\[
|\pi(A)||B\cap H|\leq |AB|.
\]
\end{proof}

We mention one more result, analogous to Petridis' version of the Pl\"{u}nnecke-Ruzsa theorem \cite{petridis2012proofs}, which will be proved in note of the author.\fxnote{add reference to proof of PPR}
\begin{prop}
  \label{prop:19}
Suppose $G\actson X$, and $A\subseteq G, Y\subseteq X$ are finite subsets.
Then there is a non-empty subset $B\subseteq A$ such that for any finite subset $C\subseteq G$, we have
\[
|CB(Y)|\leq \frac{|B(Y)||C(Y)|}{|B|}.
\]
\end{prop}

\subsubsection{Covering lemma for image sets}
Recall that Ruzsa's covering lemma states that if $|A+B|\leq K|A|$, then $B$ is covered by $K$ translates of $A-A$.
A similar statement holds image sets $|A(Y)|\leq K|Y|$, but the number of approximate orbits $A^{-1}A(x)$ needed to cover $Y$ may be much larger than $K$.
\begin{prop}[Ruzsa covering lemma for group actions]
\label{prop:2}
If \(|A(Y)|=K|Y|\), then there exists \(Z\subseteq Y\) such that
\begin{enumerate}
\item \(Y\subseteq A^{-1}A(Z)\) \label{cover}
\item \(|A(Z)|=\sum_{z\in Z}|A(z)|\) \label{disjoint}
\item \(\displaystyle |Z|\leq  \frac{K|Y|}{|A|}\left(1+\frac
    1{|A||Z|}\sum_{g\not=e} |\fix(g)\cap Z|\cdot |A\cap Ag|\right)\).\label{Z-bound}
\end{enumerate}
\end{prop}
Proposition~\ref{prop:2} recovers the Ruzsa covering lemma for multiplicative combinatorics, since the action of $G$ on itself by left-translation is \emph{free}, meaning that only the identity element $e$ has fixed points.
We defer the proof of Proposition~\ref{prop:2} to Appendix~\ref{sec:covering-lemmas}.

Write $r_{A^{-1}A}(g)=|A\cap Ag|$.
The sum in \ref{Z-bound} can be expressed as follows:
\begin{align}
  1+\frac 1{|A||Z|}\sum_{g\not=e}|\fix(g)\cap Z|r_{A^{-1}A}(g) &= \frac 1{|A||Z|}\sum_{z\in Z}\sum_{g\in\stab(z)}r_{A^{-1}A}(g)\\
& = \frac 1{|A||Z|}\sum_{z\in Z}\sum_{a\in A}|a^{-1}A\cap \stab(z)|\\
& \leq \max_{a\in A,z\in Z}|a^{-1}A\cap \stab(z)|.  
\end{align}

For a general group action, the bound for $|Z|$ in $(3)$ is worse than in the covering lemma for multiplicative combinatorics, but Example~\ref{ex:covering-sharp} in Section~\ref{sec:covering-lemmas} shows that the bound $(3)$ of Proposition~\ref{prop:2} is sharp in general.

\subsection{Symmetry sets}
\label{sec:symmetry-sets}

Let $G$ be a group acting on a set $X$.
If $g\in G$ and $Y$ is a finite subset of $X$, we say that $g$ is an \emph{$\alpha$-approximate symmetry of $Y$}\/ if $|Y\cap gY|\geq\alpha |Y|$.
The collection of all $\alpha$-approximate symmetries of a set is called a \emph{symmetry set}.
\begin{defn}[Symmetry set]
  \label{def:1}
Suppose $G\actson X$.
For $0<\alpha\leq 1$ and a finite subset $Y\subseteq X$, define the \emph{$\alpha$-symmetry set of $Y$}\/ by
\[
\sym^{G\actson X}_\alpha(Y) = \sym_\alpha(Y) = \{g\in G\colon |Y\cap g Y|\geq\alpha|Y|\}.
\]
\end{defn}
Symmetry sets were defined for abelian groups in \cite[Section 2.6]{tao2010additive}.
In other contexts, symmetry sets have been called $k$-rich transformations, where $k$ corresponds to $\alpha|Y|$; see for instance \cite{solymosi2007number,elekes2011dimension,elekes2011incidences,solymosi2012incidence,guth2015erdos}.

\begin{repex}{ex:mult-comb}[Symmetry sets in additive combinatorics]
If $G$ is an abelian group acting on itself by translation, then $\sym_\alpha^{G\actson G}(Y)$ is the set of \emph{popular differences}\/ of $Y$.
That is, \[\sym_\alpha(Y)=\{x\in G\colon r_{Y-Y}(x)\geq \alpha|Y|\},\] where $r_{Y-Y}(g)=|Y\cap (x+Y)|$ is the number of ways to write $x=y-y'$ with $y,y'\in Y$.
\end{repex}
\begin{repex}{ex:1D-aff}[Symmetry sets for affine transformations]
For $Y\subseteq\F$, the set of lines $\ell$ such that $|(Y\times Y)\cap \ell|\geq k$ is called the set of \emph{$k$-rich lines}\/ of $Y\times Y$.
If $G=\aff 1\F$ acts on $X=\F$ by affine transformations, then $\sym_\alpha^{G\actson X}(Y)$ corresponds to $k$-rich lines in $Y\times Y$ for $k=\alpha|Y|$.
That is, if $\ell$ is a non-vertical line in $\F^2$ with equation $y=gx$ for $g\in G$, then $|(Y\times Y)\cap \ell|\geq \alpha|Y|$ if and only if $|Y\cap gY|\geq\alpha|Y|$.
\end{repex}

Note that if $\alpha\geq \beta$ then $\sym_\alpha(Y)\subseteq\sym_\beta(Y)$; that is, symmetry sets of different levels are \emph{nested}.
We turn to algebraic properties of symmetry sets in the next section.

\subsubsection{Group like properties of symmetry sets}
When $\alpha=1$, the symmetry set $\sym_1(Y)$ is a subgroup of $G$: it is the stabilizer of $Y$ for the induced action of $G$ on subsets of $X$.
For $\alpha <1$, symmetry sets retain some group-like properties: $\sym_\alpha(Y)$ contains the identity and inverses ($\sym_\alpha(Y)^{-1}=\sym_\alpha(Y)$), and is approximately closed under multiplication in the following sense.
\begin{prop}[Approximate multiplicative closure]
  \label{prop:3}
If $A$ is a non-empty subset of $\sym_\alpha(Y)$, then there exists a relation $E\subseteq A^{-1}\times A$ such that
\[
|E|\geq\frac{\alpha^2}2|A|^2\andd A^{-1}\pp{E}A\subseteq\sym_{\frac{\alpha^2}2}(Y).
\]
Further, $(A^{-1}\pp{E}A)^{-1}=A^{-1}\pp{E}A$.
\end{prop}
As in the proof of the abelian case \cite[Lemma 2.33]{tao2010additive}, the proof uses Cauchy-Schwarz together with a ``popularity'' argument.

We will use the shorthand notation $Y_g:=Y\cap gY$.
The idea behind Proposition~\ref{prop:3} is that the sets $Y_g$ act, on average, like random subsets of $Y$ with density $\alpha$, so the intersections $Y_g\cap Y_{g'}$ should have density $\alpha^2$.
\begin{proof}
  Since $A\subseteq\sym_{\alpha}(Y)$, we know that
\[
\sum_{s\in A}|Y_s|\geq \alpha|Y||A|.
\]
By Lemma~\ref{lem:CSint}, there is a subset $P$ of $A\times A$ such that $|Y_s\cap Y_{s'}|\geq \frac{\alpha^2}{2}|Y|$ for all $(s,s')$ in $P$ and $|P|\geq \frac{\alpha^2}{2} |A|^2$.

Since
\[
|Y_s\cap Y_{s'}| = |Y\cap sY\cap s'Y| \leq |sY\cap s'Y| = |Y\cap s^{-1}s'Y| = |(s')^{-1}sY\cap Y|,
\]
it follows that $s^{-1}s'$ and $(s')^{-1}s$ are in $\sym_{\frac{\alpha^2}{2}}(Y)$ for all pairs $(s,s')$ in $P$.

By a slight change of notation, we may say that $s^{-1}s'\in\sym_{\alpha^2/2}(Y)$ if $(s,s')\in P$ or if $(s',s)\in P$.
Thus setting $E=P\cup P^{-1}$, we see that
\[
A^{-1}\pp E A = \{ s^{-1}s' \colon \mbox{$(s,s')\in P$ or $(s', s)\in P$}\}\subseteq\sym_{\alpha^2/2}(Y)
\]
and $(A^{-1}\pp E A)^{-1}=A^{-1}\pp E A$.
\end{proof}

The following corollary shows the utility of Proposition~\ref{prop:3}.
\begin{prop}
  \label{prop:4}
Fix $K\geq 1$ and suppose that $A\subseteq\sym_\alpha(Y)$.
If $|\sym_{\alpha^2/2}(Y)|\leq K|A|$, then there is an absolute constant $C>1$, an element $g\in A$, and a subset $S\subseteq G$ such that 
\begin{enumerate}
\item $gS\subseteq A\subseteq\sym_\alpha(Y)$ ,
\item $|S|\gg \left( \alpha/ K \right)^C|A|$,
\item $|S^3|\ll (K/\alpha)^C|S|$.
\end{enumerate}
\end{prop}
The proof of Proposition~\ref{prop:4} uses the following version of the \bsg{} theorem.
\begin{lem}
  \label{lem:1}
If $A$ and $B$ are finite subsets of a group $G$ and $E\subseteq A\times B$ is a relation such that
\[
|E|\geq \alpha|A||B| \andd |A\pp E B|\leq K|A|^{1/2}|B|^{1/2},
\]
where $\alpha\in (0,1]$ and $K>0$, then there is an element $a$ in $A$ and a subset $S\subseteq a^{-1}A$ such that
\[
|S|\gg \left( \frac \alpha K \right)^C |A|\andd |S^3|\ll \left( \frac K\alpha \right)^C|S|,
\]
where $C$ is an absolute constant.
\end{lem}
Lemma~\ref{lem:1} follows from the non-commutative \bsg{} theorem \cite[Theorem 2.44]{tao2010additive}, which yields large subsets $A'\subseteq A, B'\subseteq B$ such that $AB$ is small, and a result of Tao \cite[Proposition 4.5]{tao2008product}, which converts small doubling to small iterated growth.
\begin{proof}[Proof of Proposition~\ref{prop:4}]
Since $\sym_\alpha(Y)=\sym_\alpha(Y)^{-1}$, we have $A^{-1}\subseteq\sym_\alpha(Y)$. 
  By Proposition~\ref{prop:3} applied to $A^{-1}$, there is a subset $E\subseteq A\times A^{-1}$ such that $|E|\geq \alpha^2|A|^2/2$ and $A\pp E A^{-1}\subseteq \sym_{\alpha^2/2}(Y)$.

Thus by hypothesis we have $|A\pp E A^{-1}|\leq K|A|$ for $|E|\geq\alpha^2|A|^2/2$.
Applying Lemma~\ref{lem:1} (with $A=A, B=A^{-1}$), there is $g\in A$ and $S\subseteq G$ such that $S\subseteq g^{-1}A$, $|S|\gg (\alpha/K)^C|A|$, and $|S^3|\ll (K/\alpha)^C|A|$.
This completes the proof.
\end{proof}

Proposition~\ref{prop:4} is similar to the $\ell^2$-flattening lemma of Bourgain and Gamburd \cite{bourgain2008uniform-a}---it shows that iterating Proposition~\ref{prop:3} will yield more rich transformations unless the set of rich transformations is essentially an approximate group.

\subsubsection{Covering lemma for symmetry sets}
\label{sec:covering-lemma-symsets}

The following proposition is an approximate version of the fact that invariant sets are covering by unions of orbits.
\begin{prop}
  \label{prop:5}
If $B\subseteq\sym_\alpha(Y)$, then there exist subsets $Y'\subseteq Y$ and $Z\subseteq Y'$ such that
\begin{equation}
  \label{eq:1}
  Y'\subseteq B^{-1}B(Z),
\end{equation}
\begin{equation}
  \label{eq:2}
  |Y'| \geq \frac{\alpha|B|}{2|BB^{-1}|}|Y|,
\end{equation}
and
\begin{equation}
  \label{eq:3}
|Z|\leq \frac{2|Y|}{\alpha|B|}\max_{z\in Z}|B^{-1}B\cap\stab(z)|.
\end{equation}
In addition,
\begin{align*}
  |Z|&\leq \frac{4|Y|}{\alpha^2|B|^2}\frac 1{|Z|}\sum_{z\in Z}\sum_{g\in\stab(z)}r_{B^{-1}B}(g) \\
&= \frac{4|Y|}{\alpha^2|B|} \left( 1+\frac 1{|B||Z|}\sum_{g\not=e}|\fix(g)\cap Z|\, r_{B^{-1}B}(g) \right).
\end{align*}
\end{prop}
Even for the case of multiplicative combinatorics (Example~\ref{ex:mult-comb}), Proposition~\ref{prop:5} seems to be new, though a similar result is contained in the proof of \cite[Theorem 2.35]{tao2010additive}.
We prove Proposition~\ref{prop:5} in Appendix~\ref{sec:covering-lemmas}.

\subsubsection{Upper bounds for $|\sym_\alpha(Y)|$}
\label{sec:upper-bounds-symset}

This section contains bounds for the size of certain symmetry sets.
In conjunction with Proposition~\ref{prop:4}, these bounds imply that symmetry sets are controlled by approximate groups.
Roughly, this means that ``approximate stabilizers are approximate groups'', which is an approximate analog of the fact that stabilizers are subgroups. \fxnote{mention this ``approximate'' heuristic up front}

First, we examine the simplest bound for $|\sym_\alpha(Y)|$ in the setting of additive combinatorics.
Generalizing this bound leads to the \emph{Elekes-Sharir paradigm}, which allows us to view symmetry set bounds as a group-theoretic \emph{incidence bound}.
Then we give an example of using a product trick to get bounds 

For $G$ acting on itself by left multiplication, double counting the number of ways to write $g=y'y^{-1}$ yields
\begin{equation}
  \label{eq:4}
  |\sym_\alpha(Y)|\leq\frac{|Y|}\alpha.
\end{equation}

In general, to double count the number of ways to write $gy=y'$, we must deal with non-trivial stabilizers; $y$ and $y'$ determine $g$ only when the action is \emph{free} (as above).
To this end, we define the \emph{transporter of $y$ to $y'$}, $\trans(y,y')$, as the set of $g$ such that $gy=y'$.
The following formula is fundamental:
\begin{equation}
  \label{eq:5}
\sum_{g\in A}|Y\cap g Y| = \sum_{y,y'\in Y}|A\cap \trans(y,y')|.
\end{equation}
In particular,
\[
  |A\cap\sym_\alpha(Y)| \leq \frac 1{\alpha|Y|}\sum_{y,y'\in Y}|A\cap\trans(y,y')|.
\]
Formula (\ref{eq:5}) encapsulates the \emph{Elekes-Sharir paradigm}\/ for the Erd\H{o}s distinct distance problem \cite{elekes2011incidences,guth2015erdos}.
In \cite{elekes2011incidences}, Elekes and Sharir study the action of the 2D special Euclidean group $SE(2)$ on the real plane $\R^2$.
They note that bounding the number of rich transformations reduces to an \emph{incidence problem}\/ between ``points'' in $SE(2)$ and the ``curves'' in $SE(2)$ determined by $\trans(y,y')$.
(An equally important part of their strategy is embedding this incidence problem in $\R^3$; Guth and Katz \cite{guth2015erdos} give a parameterization where the curves $\trans(y,y')$ are straight lines in $\R^3$.)

More precisely, if $A\subseteq G$ and $P\subseteq X\times X$, we define the number of \emph{incidences}\/ between elements of $A$ and transporters $\trans(x,y)$ with $(x,y)\in P$ by
\begin{equation}
  \label{eq:elekes-sharir}
  I(P, A):=|\{((x,y),a)\in P\times A\colon a\in\trans(x,y)\}|=\sum_{(x,y)\in P}|A\cap\trans(x,y)|.
\end{equation}
Thus, (\ref{eq:5}) is an expression for $I(Y\times Y,A)$ and $\alpha|Y||A\cap\sym_\alpha(Y)|\leq I(Y\times Y,A)$ is the usual set-up for bounding the number of rich transformations.
\fxnote{references for incidence problems, esp. those related to rich transformations}

The following Proposition is a bound for the size of $\sym_\alpha(Y)$ for a general group action.
\begin{prop}
  \label{prop:6}
If $A\subseteq G$ and $Y\subseteq X$ are finite, then 
\[
  |A\cap\sym_\alpha(Y)| \leq \alpha^{-1}|Y|\max_{g\in G, x\in X}|A\cap g\stab(x)|.
\]

In particular, if $G$ is finite and $G\actson X$ transitively with point stabilizer $H$, so that $X$ is in bijection with $G/H$, then
\[
  |\sym_\alpha(Y)|\leq \alpha^{-1}|Y| \frac{|G|}{|X|}.
\]
\end{prop}
Proposition~\ref{prop:6} generalizes the bound (\ref{eq:4}) for a group acting on itself by left-translation, since for this action point stabilizers are trivial.
\begin{proof}[Proof of Proposition~\ref{prop:6}]
By (\ref{eq:5}) we have
\[
|A\cap\sym_\alpha(Y)| \leq \frac 1{\alpha|Y|} \sum_{y,y'\in Y}|A\cap\trans(y,y')|
\leq \frac{|Y|}\alpha \max_{y,y'\in Y}|A\cap \trans(y,y')|.
\]  
Since $\trans(y,y')$ is a coset of $\stab(y)$ we have
\[
|A\cap\trans(y,y')|\leq \max_{g\in G}|A\cap g\stab(y)|
\]
for all $y,y'$ in $X$, which proves the desired bound.

If $G$ is finite and $G\actson X\equiv G/H$ transitively, then $\stab(y)$ is conjugate to $H$, hence
\[
|A\cap g\stab(y)|\leq |A\cap gg'H(g')^{-1}x|\leq|H|=|G|/|X|.
\]
\end{proof}

In particular, if $G$ acts \emph{freely}\/ on $X$ (meaning that only the identity has fixed points), then $|\sym_\alpha(Y)|\leq \alpha^{-1}|Y|$.

It is often profitable to consider the diagonal action of $G$ on subsets of $X^n$ for integers $n>1$.
For sufficiently large $n$, we can sometimes obtain a free action on a subset of $X^n$.

\begin{prop}[Bound for almost-free actions]
  \label{prop:7}
Suppose that $G\actson X$ and each element of $g\not=e$ has fewer than $n$ fixed points; that is, for all $g\not=e$ in $G$, $|\fix(g)|<n$.
Then for any finite subset $Y\subseteq X$ and any $0<\alpha\leq 1$ such that $|Y|>(1+\alpha^{-1})n$, there is $0<\epsilon\leq n/\alpha(|Y|-n)$ such that
\[
|\sym_\alpha(Y)|\leq \frac{n!}{\alpha^n(1-\epsilon)}{|Y|\choose n}\leq \left( 1+\frac n{\alpha|Y|-(1+\alpha)n} \right)\frac{|Y|^n}{\alpha^n}.
\]
\end{prop}
\begin{proof}
Let $X^{(n)}\subseteq X^n$ denote the set of $n$-tuples of distinct points of $X$.

If $g\in\sym_\alpha^{G\actson X}(Y)$, then
\begin{equation}
  \label{eq:6}
  |Y^{(n)}\cap gY^{(n)}|\geq n!{\alpha|Y|\choose n}\geq \alpha^n(1-\epsilon)|Y^{(n)}|,
\end{equation}
where $\epsilon\leq n/\alpha(|Y|-n)$.
That is, $g\in\sym_{\alpha^n(1-\epsilon)}(Y^{(n)})$.

  If $|\fix(g)|<n$ for all $g\not=e$, then $G$ acts freely on $X^{(n)}$.
In particular, for any $Y\subseteq X$ and any $0<\alpha'\leq 1$ by Proposition~\ref{prop:6} we have
\[
|\sym_{\alpha'}^{G\actson X^{(n)}}(Y^{(n)})|\leq \frac{|Y^{(n)}|}{\alpha'}.
\]
By (\ref{eq:6}), we have
\[
|\sym_\alpha^{G\actson X}(Y)|\leq |\sym_{\alpha'}^{G\actson X^{(n)}}(Y^{(n)})|,
\]
where $\alpha'=\alpha^n(1-\epsilon)$.
\end{proof}
We can weaken the hypothesis on $G$ to $|\fix(g)\cap Y|<n$ for any $g\not=e$.

We end this section with a symmetry bound for $\aff 1\F\actson\F$ based on an \emph{incidence bound}\/ for points and lines in $\F^2$ (this is the standard procedure for bounding the number of ``rich lines'').

If $A\subseteq\aff 1\F$, we may associate a set of lines $L=L_A$ in $\F\times \F$ to the elements of $A$.
Let $\ell_g$ be the line corresponding to $g\in\aff 1\F$.
Then
\[
|Y\cap gY|\geq\alpha |Y|\iff |(Y\times Y)\cap \ell_g|\geq \alpha|Y|.
\]
By Cauchy-Schwarz, one can show that the number of \emph{incidences}\/ between a set of points $P\subseteq \F^2$ and a set of lines $L$ in $\F^2$ satisfies
\[
I(P,L):=|\{(p,\ell)\in P\times L\colon p\in\ell\}| \leq |P||L|^{1/2}+2|L|.
\]
Thus if $L$ denotes the set of lines corresponding to $\sym_\alpha(Y)$, we have
\[
\alpha|Y||\sym_\alpha(Y)| \leq I(Y\times Y,L)\leq |Y|^2|\sym_\alpha(Y)|^{1/2}+2|\sym_\alpha(Y)|.
\]
Thus if $\alpha|Y|>2$, we have
\[
|\sym_\alpha(Y)|\leq \frac{|Y|^2}{(\alpha-2/|Y|)^2}.
\]
This is just the Elekes-Sharir paradigm (\ref{eq:5}): if $y,y'\in \F$, then $\trans(y,y')$ corresponds to a pencil of lines in $\F^2$, and $g\in\trans(y,y')$ if $g$ is contained in this pencil; by duality, this corresponds to point-line incidences.

A similar argument works whenever we have an incidence bound that corresponds to incidence relation $g\in\trans(y,y')$.

\subsection{Action Energy}
\label{sec:energy}

In this section, we generalize \emph{multiplicative energy}\/ to group actions.
Recall that if $A$ and $B$ are finite subsets of a group $G$, then the multiplicative energy of $A$ and $B$ is the number of multiplicative quadruples $ab=a'b'$ with $a,a'\in A, b,b'\in B$:
\[
E(A,B):=|\{(a,a',b,b')\in A\times A\times B\times B\colon ab=a'b'\}|.
\]

\begin{defn}[Action energy]
\label{def:2}
If $G\actson X$ and $A\subseteq G, Y\subseteq X$ are finite, then the \emph{action energy} of \(A\) and \(Y\) is defined by
\[
E(A,Y):=|\{(a_1,a_2,y_1,y_2)\in A\times A\times Y\times Y\colon a_1y_1=a_2y_2\}|.
\]
\end{defn}
For $G$ acting on itself by left translation, the action energy is multiplicative energy (or additive energy if $G$ is abelian and written additively).
Another instance of action energy occurs in \cite{yazici2015growth}: given a set of affine transformations $A$ in $\aff 1\F$ and a subset $Y\subseteq\F$, the number of \emph{collisions}\/ of image of lines was defined as the number of solutions to
\[
a_1(y_1)=a_2(y_2)
\]
with $a_1,a_2\in A$ and $y_1,y_2\in Y$.
This is the action energy $E(A,Y)$ for $G=\aff 1\F$ acting on $X=\F$ by affine transformations.\fxnote{use example environment for energy?}

\subsubsection{Alternate expressions for action energy}
As with additive and multiplicative energy, there are many useful expressions for $E(A,Y)$.
We may express \(E(A,Y)\) as a sum:
 \begin{equation}
   \label{eq:7}
   E(A,Y)=\sum_{a_1,a_2\in A}|a_1 Y\cap a_2 Y|.
 \end{equation}
Setting
\[
\rep{A}(g) = |\{(a_1,a_2)\in A\times A\colon a_1^{-1}a_2 = g\}|,
\]
we have
\begin{equation}
  \label{eq:8}
  E(A,Y) = \sum_{g\in G}\rep{A} (g)|Y\cap gY|.
\end{equation}
We can further decompose (\ref{eq:8}) using transporters:
\begin{equation}
  \label{eq:9}
  E(A,Y)=\sum_{g\in G} \rep{A}|Y\cap gY|=\sum_{y,y'\in Y}\sum_{g\in\trans(y,y')}r_{A^{-1}A}(g).
\end{equation}

Setting
\[
r_{A(Y)}(x)=|\{(a,y)\in A\times Y\colon a(y)=x\}|,
\]
we have
\begin{equation}
  \label{eq:10}
  E(A,Y)=\sum_{x\in A(Y)} r_{A(Y)}^2(x).
\end{equation}
Note that if $\stab(x)=\{e\}$, then $r_{A(Y)}(x)=|Y\cap A^{-1}(x)|$; in general,
\[
r_{A(Y)}(x) = \sum_{y\in Y}|A\cap\trans(y,x)| = \sum_{y\in Y\cap A^{-1}(x)}|A\cap a_y\stab(x)|,
\]
where for each $y$ in $Y\cap A^{-1}(x)$, we have chosen an element $a_y$ of $A$ such that $x=a_y(y)$; that is, $a_y\in A\cap\trans(y,x)$.

\subsubsection{Bounds for $E(A,Y)$}

In this subsection we record some upper bounds for $E(A,Y)$.
For $G$ acting on itself by left translation, we have
\begin{equation}
  \label{eq:11}
  E(B,A) \leq |B|^2|A|
\end{equation}
and
\begin{equation}
  \label{eq:12}
  E(B,A) \leq |B||A|^2
\end{equation}
since any three variables in the equation \[b_1a_1=b_2a_2\] determines the fourth.

For a general $G$-set $X$, given $g_1,g_2$ in $G$ and $x_1$ in $X$, the equation \[g_1x_1=g_2x_2\] determines $x_2$, thus
\begin{equation}
  \label{eq:13}
  E(A,Y)\leq |A|^2|Y|.
\end{equation}
However, if we are given $x_1$ and $x_2$, then all we know is that $g_2^{-1}g_1x_1=x_2$, which means that $g_2^{-1}g_1$ is in the transporter of $x_1$ to $x_2$; if $|\trans(x_1,x_2)| > 1$, then $g_2^{-1}g_1$ is not determined by $x_1$ and $x_2$.
Thus, \eqref{eq:11} generalizes, but \eqref{eq:12} does not.

The following proposition shows that the energy $E(A,Y)$ is large if and only if many elements of $A^{-1}A$ are contained in a symmetry set of $Y$.
\begin{prop}[Generic upper and lower bounds for action energy]
\label{prop:8}
For any $\alpha$ in $(0,1]$, we have
\begin{equation}
  \label{eq:14}
  \alpha|Y|\sum_{g\in \sym_\alpha(Y)}r_{A^{-1}A}(g) \leq E(A,Y)
\end{equation}
and
\begin{align}
    E(A,Y)&\leq (\ceil{\alpha |Y|}-1)|A|^2+|Y|
\sum_{g\in \sym_\alpha(Y)}r_{A^{-1}A}(g) \label{eq:15}\\  
&\leq (\ceil{\alpha |Y|}-1)|A|^2+|A||Y| \max_{a\in A}|a^{-1}A\cap\sym_\alpha(Y)|. \label{eq:47b}
\end{align}
\end{prop}
\begin{proof}
  To prove \eqref{eq:14}, we use expression~\eqref{eq:8} for $E(A,Y)$ and restrict the sum to transformations $g$ in $\sym_\alpha(Y)$.

To prove \eqref{eq:15}, we split the sum in expression~\eqref{eq:8} into two terms:
\[
E(A,Y) = \sum_{g\not\in\sym_{\alpha}(Y)}\rep{A}(g)|Y\cap gY|+ \sum_{g\in\sym_{\alpha}(Y)}\rep{A}(g)|Y\cap gY| = I + II.
\]
To bound sum $I$, we use the upper bound $|Y\cap gY| < \ceil{\alpha |Y|}$ and sum over $g$, and to bound sum $II$ we use the bound $|Y\cap gY|\leq |Y|$.

To prove \eqref{eq:47b}, write $S=\sym_\alpha(Y)$.
Then
\[
\sum_{g\in \sym_\alpha(Y)}r_{A^{-1}A}(g)= \sum_{g\in S}r_{A^{-1}A}(g)=\sum_{a\in A}r_{AS}(a)=\sum_{a\in A}|a^{-1}A\cap S|,
\]
and the bound follows by pigeonholing over $a$.
\end{proof}

The generic upper bound implies an approximate version of the orbit-stabilizer theorem, which is frequently used in approximate group theory (see for example, \cite{helfgott2015growth}).
\begin{cor}[Orbit-stabilizer theorem for sets]
\label{cor:1}
 Let \(G\) be a group acting on a set \(X\).
Fix \(x\) in \(X\) and let \(A\) be a non-empty subset of \(G\).
Then there is an element $a_0$ in $A$ such that
\[
|a_0^{-1}A\cap\stab(x)||A(x)|\geq |A|.
\]
In fact,
\[
\frac{|A|^2}{|A(x)|}\leq \sum_{g\in\stab(x)}r_{A^{-1}A}(g) = \sum_{a\in A}|a^{-1}A\cap\stab(x)|.
\]
\end{cor}
\begin{proof}
Apply Proposition~\ref{prop:8} with $Y=\{x\}$ and $\alpha =1$
Then $\sym_\alpha(Y) = \stab(x)$ and $\ceil{\alpha|Y|}-1=0$, so we have
\[
E(A, \{x\})\leq \sum_{g\in\stab(x)}r_{A^{-1}A}(g)
\leq |A| \max_{a\in A}|a^{-1}A\cap \stab(x)|.
\]
By equation~\eqref{eq:21} (below), we have
\[
|A|^2\leq |A(x)|\cdot E(A,\{x\}),
\]
thus the desired result follows by combining these two inequalities and choosing $a_0$ to be argument of the maximum.
\end{proof}
In applications, $A$ is typically an approximate group, so the number of solutions to $a_2^{-1}a_1=g$ is near maximal for all $g$ in $A^{-1}A$, and this inequality is not as wasteful as it may seem.

\fxnote{orphan remark on non-trivial bound for number of collisions}
We end with a remark on the existence of \fxnote*{what do we mean by non-trivial?}{non-trivial} upper bounds for $E(A,Y)$.
For the energy of the 1D affine group over $\F_p$ acting on $\F_p$, if $|A||Y|\leq p^2$, then it was proved in \cite{yazici2015growth} that
\[
E(A,Y) \ll (|A||Y|)^{3/2}+\max(\kappa, |A|)|A||Y|,
\]
where $\kappa$ is the maximum number of elements in $A$ that are contained in a coset of an abelian subgroup of the affine group.
\fxnote{give asymptotic formula for collisions bound and drop constraints}
\subsection{Conversions between image sets, symmetry sets, and action energy}
\label{sec:conversions}

In additive combinatorics, there are many ways of quantifying the ``additive structure'' of a set: small sumset, large additive energy, many popular differences.
It is possible to convert between many of these forms.
In this section, we give the analogous conversions for image sets, symmetry sets, and action energy.

We show that the conditions
\begin{align}
&E(A,Y)\geq \alpha|A|^2|Y|\label{eq:16}\\
 &A\subseteq\sym_\alpha(Y)\label{eq:17}\\
 &|A_E(Y)|\leq K|Y|, \quad\mbox{where $|E|\geq\rho|A||Y|$}\label{eq:18}
\end{align}
are equivalent with polynomial dependence on parameters (for example, if (\ref{eq:18}) holds, then (\ref{eq:16}) and (\ref{eq:17}) are true with $\alpha$ polynomial in $K$ and $\rho$).

We will show that (\ref{eq:17}) and (\ref{eq:18}) imply (\ref{eq:16}) and conversely, which shows indirectly that (\ref{eq:17}) and (\ref{eq:18}) are equivalent.
At the end of the section, we prove directly that (\ref{eq:17}) and (\ref{eq:18}) are equivalent.

To prove that (\ref{eq:18}) implies (\ref{eq:16}), we use a Cauchy-Schwarz lower bound for $E(A,Y)$ analogous to the Cauchy-Schwarz lower bound for multiplicative energy
\[
|B|^2|A|^2\leq |BA|\cdot E^\times(B,A),
\]
which says that the size of the product set and multiplicative energy are inversely correlated.

To prove the corresponding bound for action energy, we introduce some notation.
If $E\subseteq A\times Y$, we define
\[
r_E(x):=|\{(a,y)\in E\colon a(y)=x\}|.
\]
We have
\begin{equation}
  \label{eq:19}
  |A||Y|=\sum_x r_{A(Y)}(x)
\end{equation}
and
\begin{equation}
  \label{eq:20}
  |E|=\sum_x r_E(x).
\end{equation}
The Cauchy-Schwarz bound for action energy is the following.
\begin{prop}[Small image set implies large energy]
\label{prop:9b}
For any subset \(A\) of \(G\) and any subset \(Y\) of \(X\), we have
\begin{equation}
  \label{eq:21}
  |A|^2|Y|^2\leq |A(Y)|E(A,Y).
\end{equation}
Further, for any subset \(E\subseteq A\times Y\), we have
\begin{equation}
  \label{eq:22}
  | E|^2 \leq |A_E(Y)|E(A,Y).
\end{equation}
\end{prop}
The proof is a single application of the Cauchy-Schwarz inequality to (\ref{eq:21}) or (\ref{eq:22}), using $r_{E}(x)\leq r_{A(Y)}(x)$ and expression~(\ref{eq:10}) for $E(A,Y)$.

By Proposition~\ref{prop:9b}, condition (\ref{eq:18}) implies that
\begin{equation}
  \label{eq:23}
  E(A,Y)\geq\frac{\rho^2|A|^2|Y|}{K},
\end{equation}

Now we prove that (\ref{eq:17}) implies (\ref{eq:16}) by comparing $E(A,Y)$ to $\sym_\alpha(Y)$:
\begin{prop}
  \label{prop:9}
If $G\actson X$ and $A\subseteq G, Y\subseteq X$ are finite, then for all $\alpha\in [0,1]$
\[
  \alpha^2 |A\cap \sym_\alpha(Y)|^2 |A| \leq E(A,Y).
\]
\end{prop}
\begin{proof}
Let $S=A\cap \sym_\alpha(Y)$.
We have
\[
\alpha|S||Y|\leq \sum_{g\in S}|Y_g|,
\]  
thus by the Cauchy-Schwarz intersection lemma (Lemma~\ref{lem:CSint}) and expression (\ref{eq:7}) for $E(A,Y)$,
\[
\alpha^2|S|^2|Y| \leq \sum_{g,g'\in S}|Y_g\cap Y_{g'}|\leq \sum_{g,g'\in A}|gY\cap g'Y| = E(A,Y).
\]
\end{proof}
By Proposition~\ref{prop:9}, condition (\ref{eq:17}) implies that
\begin{equation}
  \label{eq:24}
  E(A,Y)\geq \alpha^2|A|^2|Y|.
\end{equation}

Now we will prove the reverse implications.
\begin{prop}
  \label{prop:10}
Suppose that $E(A,Y)\geq 2\alpha|A|^2|Y|$.
Then
\begin{enumerate}
\item there exists $E\subseteq A\times Y$ such that $|E|\geq \alpha|A||Y|$ and  $|A_E(Y)|\leq\alpha^{-2}|Y|$,
\item there exists $a_0$ in $A$ such that $|a_0^{-1}A\cap\sym_\alpha(Y)|\geq\alpha|A|$.
\end{enumerate}
\end{prop}
\begin{proof}
To prove 2, we use Proposition~\ref{prop:8}:
\[
2\alpha|A|^2|Y| \leq E(A,Y) \leq \alpha|A|^2|Y| + |Y||A|\max_{a\in A}|a^{-1}A\cap\sym_\alpha(Y)|,
\]
hence there is an $a_0$ in $A$ such that
\[
\alpha|A| \leq |a_0^{-1}A\cap\sym_\alpha(Y)|.
\]

To prove 1, we use a popularity argument.
Let $P\subseteq A(Y)$ denote the set of $x$ such that $r_{A(Y)}(x)\geq \alpha|A|$.
Since
\[
2\alpha|A|^2|Y| \leq E(A,Y) \leq \sum_{x\in P}r_{A(Y)}^2(x) + \alpha|A| \sum_{x\not\in P}r_{A(Y)}(x),
\]
we have
\[
\alpha|A|^2|Y| \leq \sum_{x\in P}r_{A(Y)}^2(x).
\]

Now, let $E$ denote the set of pairs $(a,y)\in A\times Y$ such that $a(y)\in P$.
Since $r_{A(Y)}(x)\leq |A|$, we have
\[
|E| = \sum_{x\in P}r_{A(Y)}(x) \geq \alpha|A||Y|.
\]
By definition, $A_E(Y)\subseteq P$, so to bound $|A_E(Y)|$ it suffices to bound $|P|$.

Since 
\[
\alpha^2|A|^2|P| \leq |A|^2|Y|,
\]
we have $|P|\leq\alpha^{-2}|Y|$, which proves the desired bound on the partial image set.
\end{proof}

It is possible to prove the equivalence of (\ref{eq:17}) and (\ref{eq:18}) directly.
We briefly sketch proofs.
\begin{itemize}
\item Suppose that $|A_E(Y)|\leq K|Y|$ and $|E|\geq \rho|A||Y|$.
Let $A'$ denote the set of $a\in A$ such that $(a,y)\in E$ for at least $\rho|Y|/2$ elements $y\in Y$; by popularity, we have $|A'|\geq \rho|A|/2$.
Let $Y':=Y\cup A_E(Y)$; by assumption $|Y'|\leq(K+1)|Y|$.
On the other hand, for each $a\in A'$ we have \[|Y'\cap aY'|\geq \rho|Y|/2\geq \rho|Y'|/2(K+1).\]
Thus for $\alpha=\rho/2(K+1)$ we have $A'\subseteq\sym_\alpha(Y')$, with $|Y|\leq |Y'|\leq (K+1)|Y|$ and $|A'|\geq \rho|A|/2$.
\item Conversely, if $A\subseteq\sym_\alpha(Y)$, then the set $E$ of pairs $(a,y)$ in $A\times Y$ such that $a(y)\in Y$ has size at least $\alpha|A||Y|$, and $A_E(Y)\subseteq Y$.
\end{itemize}

\addtocontents{toc}{\protect\setcounter{tocdepth}{1}}
\section{Group action \bsg{} theorem}
\label{sec:group-action-bsg}

Let $A$ and $B$ be finite subsets of an additive group.
By Cauchy-Schwarz, if the sum set $A+B$ is small, then the additive energy $E^+(A,B)$ is large.
The \emph{\bsg{} theorem}\/ provides a partial converse in the \emph{symmetric case}, where $|A|$ and $|B|$ are roughly equal.
\begin{thm}[Symmetric \bsg{} theorem]
  \label{thm:bsg-additive-balanced}
  If $A$ and $B$ are finite subsets of an additive group and $E^+(A,B)\geq |A|^{3/2}|B|^{3/2}/K$, then there are subsets $A'\subseteq A$ and $B'\subseteq B$ such that $|A'|\gg K^{-C}|A|, |B'|\gg K^{-C}|B|$, and $|A'+B'|\ll K^C|A|^{1/2}|B|^{1/2}$, where $C>1$ is an absolute constant.
\end{thm}
The upper bound $E^+(A,B)\leq |A|^2|B|$ implies that the hypotheses of Theorem~\ref{thm:bsg-additive-balanced} hold only if $K\geq |B|^{1/2}/|A|^{1/2}$, which makes the upper bound on $|A'+B'|$ is trivial if $A$ is much smaller than $B$.
In this case, Tao and Vu proved a non-trivial converse to the
Cauchy-Schwarz lower bound for additive energy, known as the
\emph{asymmetric \bsg{} theorem} \cite[Theorem 2.35]{tao2010additive}.
\begin{thm}[Asymmetric \bsg{} theorem]
  \label{thm:2}
Let $A$ and $B$ be finite subsets of an additive group $Z$ such that $E(A,B)\geq 2\alpha |A|^2|B|$ and $|B|\leq L|A|$ for some $L\geq 1$ and $0<\alpha\leq 1$.

For all $\epsilon>0$ there is a constant $C_\epsilon>0$ such that the following holds:
\begin{enumerate}
\item there is a set $H$ in a translate of $\sym_{\alpha^{C_\epsilon}}(B)$ such that $|H|\leq C_\epsilon\alpha^{-C_\epsilon}|B|$ and $|H+H|\leq C_\epsilon\alpha^{-C_\epsilon}L^\epsilon|H|$,
\item there is an element $x$ in $Z$ such that $|A\cap (x+H)|\geq C_\epsilon \alpha^{C_\epsilon}L^{-\epsilon}|A|$
\item there is a subset $X$ of $Z$ with cardinality $|X|\leq C_\epsilon \alpha^{-C_\epsilon}L^\epsilon |B|/|H|$ such that $|B\cap (X+H)| \geq C_\epsilon \alpha^{C_\epsilon}L^{-\epsilon}|B|$.
\end{enumerate}
\end{thm}
Part (1) says that some translate of a symmetry set of $B$ contains a set $H$ with small doubling, part (2) says that a large portion of $A$ is contained in $H$, and Part (3) says that a large part of $B$ is covered by translates of $H$.

In the symmetric/balanced case, the \bsg{} theorem only loses powers of $\alpha$; in the asymmetric case, where $A$ is much smaller than $B$, we also lose powers of $L$.
If $L\approx 1$, then taking $A' = A\cap (x+H)$ and $B'=B\cap(x_0+H)$ for some $x_0$ in $X$, we have $|A'|,|B'|\gg_\alpha \min(|A|,|B|)$ and $|A'B'|\ll_\alpha \min(|A|,|B|)$.
However, if $L\ggg 1$, then we only have $|B\cap (x_0+H)|\gg_\epsilon \alpha^{2C_\epsilon}L^{-2\epsilon}|H|$ and $|A'+B'|\ll_\epsilon \alpha^{-C_\epsilon}L^\epsilon|B'|$.

In this section, we prove a version of the asymmetric \bsg{} theorem for \emph{group actions}, with the translation action of an abelian group $Z$ on itself replaced by a general action of a group $G$ on a set $X$.

Instead of considering sets $A,B\subseteq Z$ such that $E^+(A,B)\geq 2\alpha|A|^2|B|$, we will consider sets $A\subseteq G$ and $Y\subseteq X$ such that $E(A,Y)\geq 2\alpha|A|^2|Y|$, or equivalently, such that $A\subseteq\sym^{G\actson X}_\alpha(Y)$.

In Theorem~\ref{thm:3}, we introduce an integer parameter $J\geq 0$ such that $1/J$ plays the same role as $\epsilon$ in Theorem~\ref{thm:2}.
Theorem~\ref{thm:3} has three parts, corresponding to the three parts of Theorem~\ref{thm:2}:
\begin{itemize}
\item Part (1) says that there is a set $H$ in $\sym^{G\actson X}_{\alpha_J}(Y)$, where $\alpha_J=2(\alpha/2)^{2^J}$, with small tripling; in contrast to Theorem~\ref{thm:2}, we have no cardinality bound for $|H|$.
\item Part (2) is exactly the same: $A$ has large overlap with a translate of $H$.
\item Part (3) is similar: a large part of $Y$ is covered by images of $H$, but we do not as strong a bound for the number of images $Y$ as in Theorem~\ref{thm:2}. 
\end{itemize}
In parts (2) and (3), we introduce some additional flexibility: we may replace $H$ by a dense subset $S\subseteq H$.
We do this so that if $|H^3|\leq K|H|$ implies that $H$ has large overlap with some ``structured set'' $S$, then this structure can be transferred to $A$ and $Y$.

Recall that $A_{(n)}=(A\cup A^{-1}\cup\{e\})^n$.

\begin{thm}[asymmetric \bsg{} theorem for group actions]
\label{thm:3}
There is an absolute constant $C>0$ such that the following holds.

Let $G$ be a group acting on a set $X$, let $Y$ be a finite subset of $X$, and let $A$ be a finite subset of $G$.

Given a number $0<\alpha<1$ and an integer $J\geq 0$, let $\alpha_J=2(\alpha/2)^{2^J}$, and define $K>0$ by
\begin{equation}
  \label{eq:25}
  K = \left( \frac{|\sym_{\alpha_J}(Y)|}{|A|} \right)^{1/J}.
\end{equation}

If $A\subseteq\sym_\alpha(Y)$, then
\begin{enumerate}
\item there is an element $g_*$ in $A_{(2^J)}$ and a finite subset $A_*\subseteq G$ such that
\begin{equation}
  \label{eq:26}
g_*^{-1}  A_*\subseteq A_{(2^J)}\cap\sym_{\alpha_J}(Y)
\end{equation}
and
\begin{equation}
  \label{eq:27}
  |A_*^3| \ll \left( \frac{K}{\alpha_J}  \right)^C|A_*|,
\end{equation}
\item for any subset $S\subseteq G$ there is an element $g$ in $A_{(2^{J+1})}$ such that
\begin{equation}
  \label{eq:28}
  |A^{-1}\cap gS| + |A\cap gS| \gg \left( \frac{\alpha_J}{K}  \right)^C\frac{|S\cap A_*|}{|A_*|} |A|.
\end{equation}
\item for any $S\subseteq G$, if $B=A_*\cap S$ and $\rho=|B|/|A_*|$, then there is  a subset $Y'\subseteq Y$ such that
\begin{equation}
  \label{eq:29}
  |Y'|\geq \frac{\alpha_J|B|}{2|BB^{-1}|}|Y|\gg \left( \frac{\rho\cdot\alpha_J}{K}  \right)^C |Y|,
\end{equation}
and a subset $Z\subseteq Y'$ such that $Y'\subseteq B^{-1}B(Z)$ and
\begin{align}
  \label{eq:30}
|Z| &\leq \frac{2|Y|}{\alpha_J|B|} \left( 1+\frac 1{|B||Z|}\sum_{g\not=e}|\fix(g)\cap Z| r_{B^{-1}B}(g) \right)\\
&\ll \left( \frac{K}{\rho\cdot\alpha_J}  \right)^C  \frac{|Y|}{|S|} \left( 1 + \frac 1{|S||Z|}\sum_{g\not=e}|\fix(g)\cap Z|r_{S^{-1}S}(g) \right). \nonumber
\end{align}
\fxnote{since $S$ is an approximate group, perhaps including $r_{S^{-1}S}(g)$ is unnecessary? ...need $B$ instead of $S$... might be useful since Larson-Pink applies...}
\end{enumerate}
\end{thm}

\subsection{Preliminaries}
Now we state the lemmas necessary for the proof of Theorem~\ref{thm:3}.
The first is a uniform version of the approximate closure property of symmetry sets (Proposition~\ref{prop:3}).

Recall that if $A$ and $B$ are finite subsets of a group and $E\subseteq A\times B$ then
\[
A\pp{E} B :=\{ab\colon (a,b)\in E\}\andd r_E(x):=|\{(a,b)\in E\colon ab=x\}|.
\]
\begin{lem}[Uniform approximate closure]
  \label{lem:2}
If $A$ is a non-empty subset of $\sym_\alpha(Y)$ then there is a relation $E\subseteq A^{-1}\times A$ such that
\begin{align}
  |E| &\geq \frac{\alpha^2}{2+2\log(|A|)}\,|A|^2,\label{eq:33}\\
r_E(x) &\geq\frac{|E|}{2|A^{-1}\pp{E}A|}\qquad\mbox{for all $x$ in $A^{-1}\pp{E}A$},\label{eq:34}\\
& A^{-1}\pp{E} A \subseteq \sym_{\frac{\alpha^2}2}(Y)\label{eq:35}.
\end{align}
Further, $E$ is symmetric, so that 
\begin{equation}
  \label{eq:36}
  (A^{-1}\pp E A)^{-1}=A^{-1}\pp E A.
\end{equation}
\end{lem}
The proof of Lemma~\ref{lem:2} is essentially the same as the proof of \cite[Lemma 2.34]{tao2010additive}: combine Proposition~\ref{prop:3} with a dyadic pigeonholing argument.
(Although \cite[Lemma 2.34]{tao2010additive} is stated for abelian groups, the proof works verbatim for non-abelian groups.)
\begin{proof}[Proof of Lemma~\ref{lem:2}]
  By Proposition~\ref{prop:3}, there exists $E\subseteq A^{-1}\times A$ such that $E=E^{-1}$ (that is, $(x,y)\in E$ if and only if $(y,x)\in E$), $|E|\geq\alpha^2|A|^2/2$, and $A^{-1}\pp E A\subseteq\sym_{\alpha^2/2}(Y)$.

Let $B_j=\{x\colon 2^j\leq r_E(x) < 2^{j+1}\}$ and let $E_j=\{(a^{-1},a')\in E\colon a^{-1}a'\in B_j\}$.
Since $r_E(x)\leq |A|$ for all $x$, we have
\[
|E| = \sum_{j=0}^{\log_2|A|}\sum_{x\in B_j}r_E(x).
\]
By pigeonholing, there is a $j$ such that
\[
|E_j|=\sum_{x\in B_j}r_E(x) \geq \frac{|E|}{1+\log_2 |A|}.
\]

Since $E=E^{-1}$, we have $r_E(x)=r_E(x^{-1})$, hence $B_j=B_j^{-1}$. 
Further, $A^{-1}\pp{E_j}A=B_j$.
By definition of $B_j$,
\[
2^j|B_j|\leq \sum_{x\in B_j}r_E(x) = |E_j| < 2^{j+1}|B_j|,
\]
so for all $x$ in $B_j$,
\[
r_E(x) \geq 2^j > \frac{|E_j|}{2|B_j|}.
\]
Setting $E=E_j$ completes the proof.
\end{proof}

Lemma~\ref{lem:2} implies that if a set $S$ is dense in the product set $A^{-1}\pp E A$, then some translate of $S$ is dense in $A$.
Thus, if we find a ``structured'' subset of the product set $A^{-1}\pp E A$, we may bring that structure back to the original set $A$.
\begin{lem}[Bringing structure back]
  \label{lem:3}
If $A$ is a finite subset of $G$ and $E\subseteq A^{-1}\times A$ satisfies (\ref{eq:33}) and (\ref{eq:34}), then for any subset $S$ of $G$, there is an element $a$ in $A$ such that
\[
\frac{|A\cap aS|}{|A|}\geq \frac{\alpha^2}{4(1+\log(\alpha^{-1}|A|))}\frac{|(A^{-1}\pp E A)\cap S|}{|A^{-1}\pp E A|}.
\]
\end{lem}
\begin{proof}
Count the number of solutions to $a'=as$:
\[
\sum_{a\in A}|A\cap aS|
\geq \sum_{s\in S}r_E(s)
\geq \frac{|(A^{-1}\pp E A)\cap S| |E|}{2|A^{-1}\pp E A|}.
\]
Pigeonhole over $a$ and use (\ref{eq:33}) to complete the proof.
\end{proof}

\subsection{Proof of Theorem~\ref{thm:3}}
\label{sec:proof-group-action-bsg}

We begin by fixing notation and defining a sequence of sets inductively using Lemma~\ref{lem:2} and Lemma~\ref{lem:3}.

Given a number $0<\alpha\leq 1$ and an integer $J\geq 0$, define $\alpha_j = 2(\alpha/2)^{2^j}$ so that $\alpha_0=\alpha$ and $\alpha_{j+1}=\alpha_j^2/2$ for $j=0,\ldots,J-1$.
Define a sequence of sets $A_j\subseteq\sym_{\alpha_j}(Y)$ inductively by setting $A_0=A\cup A^{-1}\cup\{e\}$ and setting
\[
A_{j+1}:= A_j^{-1}\pp{E_j} A_j
\]
for $j=0,\ldots, J-1$, 
where $E_j\subseteq A_j^{-1}\times A_j$ satisfies (\ref{eq:33}), (\ref{eq:34}), and (\ref{eq:35}).
Such an $A_{j+1}$ exists by Lemma~\ref{lem:2}.
Further, $A_j^{-1}=A_j$ for $j=0,\ldots, J$ by \eqref{eq:36}.

Define \[L_j = \frac 1{1+\log |A_j|}.\]
For $j=0,\ldots, J$, we have
\begin{equation}
  \label{eq:37}
  A_j \subseteq A_{(2^j)}\cap\sym_{\alpha_j}(Y),
\end{equation}
which gives us the rudimentary bound $|A_j|\leq |A_0|^{2^j}$.
Hence
\begin{equation}
  \label{eq:38}
 L_j \geq 2^{-j}L_0 \andd |E_j| \geq 2^{-j}\alpha_{j+1}L_0|A_j|^2.
\end{equation}
 
By Lemma~\ref{lem:1}, for any subset $S\subseteq G$ there is an element $a_j$ in $A_j$ such that
\begin{equation}
  \label{eq:39}
  \frac{|A_j\cap a_j S|}{|A_j|} \geq \frac{\alpha_{j+1}L_0}{2^j} \frac{|A_{j+1}\cap S|}{|A_{j+1}|}.
\end{equation}

\begin{proof}[Proof of Part (1)]
Now we will show that there are two consecutive terms $A_j, A_{j+1}$ of the sequence with comparable size, hence $A_j$ will have small partial doubling.
Define $K>0$ by
\[
K^J:=\frac{|A_J|}{|A_0|}.
\]
By the pigeonhole principle, there is an index $0\leq j < J$ such that
  \begin{equation}
    \label{eq:40}
    |A_j^{-1}\pp{E_j}A_j| \leq K|A_j|,
  \end{equation}
since
\[
\prod_{j=0}^{J-1}\frac{|A_j^{-1}\pp{E_j}A_j|}{|A_j|} = \prod_{j=0}^{J-1}\frac{|A_{j+1}|}{|A_j|} = \frac{|A_J|}{|A_0|} = K^J.
\]

Now we convert from small partial doubling to small tripling.
By Lemma~\ref{lem:1} with $A=A_j^{-1},B=A_j,E=E_j, \alpha\gg\alpha_{j}^2L_j$, and $K$ as above, it follows that
there is an element $g_*\in A_j$ and a subset $A_*\subseteq g_*A_j^{-1}=g_*A_j$ such that
\begin{equation}
  \label{eq:41}
  |A_*^3| \ll \left( \frac K{L_j\alpha_j} \right)^C|A_*|
\end{equation}
and
\begin{equation}
  \label{eq:42}
  |A_*| \gg  \left( \frac{L_j\alpha_j}K \right)^C|A_{j}|,
\end{equation}
where $C>0$ is an absolute constant.
This proves part (1) of Theorem~\ref{thm:3}.
\end{proof}

\begin{proof}[Proof of Part (2)]
Now we prove part (2) of Theorem~\ref{thm:3}.
First, we will show that for any $S\subseteq G$ and any $i\in\{0,\ldots,J\}$, there exists an element $g$ in $A_{(2^i)}$ such that
\begin{equation}
  \label{eq:43}
  |A_0\cap gS|\gg \alpha_{i}^2L_i\frac{|A_i\cap S|}{|A_i|} |A_0|.
\end{equation}
This implies (2), since by (\ref{eq:42})
\[
\frac{|A_j\cap g_*^{-1}(A_*\cap S)|}{|A_j|} =\frac{|A_*\cap S|}{|A_j|}  \gg \left( \frac{L_j\alpha_{j}}K \right)^{C}\frac{|A_*\cap S|}{|A_*|},
\]
and $g_*^{-1}g\in A_{(2^j)}^2\subseteq A_{(2^{j+1})}\subseteq A_{(2^{J+1})}$.


To prove (\ref{eq:43}) we use (\ref{eq:39}) and set $g=a_0\cdots
a_{j-1}$: 
\[
|A_0\cap gS| \gg \left( \prod_{j=0}^{J-1} \alpha^2_jL_j\right)\frac{|A_j\cap S|}{|A_j|} |A_0|.
\]
This proves the claim, since
\[
\prod_{j=0}^{i-1} \alpha^2_jL_j \gg L_0 \prod_{j=0}^{i-1} 2^{-j}\left( \frac{\alpha^2_0}2\right)^{2^j} \gg L_0 2^{-i}\alpha_i^2
\]
and by (\ref{eq:37})
\[
g = a_0\cdots a_{i-1} \in (A\cup A^{-1})\cdots (A\cup A^{-1})^{2^{i-1}} \subseteq (A\cup A^{-1})^{2^i}.
\]
\end{proof}

\begin{proof}[Proof of Part (3)]
Recall that $S\subseteq G$, $B:=A_*\cap S$, and $\rho:=|B|/|A_*|$.

Since $B\subseteq A_*$, we have
\begin{equation}
  \label{eq:44}
  g_*^{-1}B\subseteq A_{(2^J)}\cap\sym_{\alpha_J}(Y),
\end{equation}
and
\begin{equation}
  \label{eq:45}
  |B^3|\ll (\alpha_J^{-1}K)^{C}|A_*| \ll  (\alpha_J^{-1}K)^{C}\rho^{-1}|B|.
\end{equation}
Applying Proposition~\ref{prop:5} to $B_0=g_*^{-1}B$, we have subsets $Y'\subseteq Y$ and $Z\subseteq Y'$ such that $Y'\subseteq B^{-1}B(Z)$,
\[
|Y'| \geq \frac{\alpha_J|B_0|}{2|B_0B_0^{-1}|}|Y|
\]
and
\[
|Z| \leq \frac {2}{\alpha_J}\frac{|Y|}{|B_0|} \left(1+ \frac 1{|B_0||Z|}\sum_{g\not=e}|\fix(g)\cap Z| r_{B_0^{-1}B_0}(g) \right).
\]
Since $|B_0B_0^{-1}|=|BB^{-1}|\ll (K/\alpha_J)^C|B|$, the desired conclusion follows.
\end{proof}

Finally, note that
\[
\frac 1{L_j} \ll 2^j \log(|A|) \ll 2^J K^\epsilon
\]
for all $\epsilon>0$, provided that $|A|\gg_\epsilon 1$ is sufficiently large; thus we may ignore the term $L_j$ by increasing the constant $C$ slightly.
\addtocontents{toc}{\protect\setcounter{tocdepth}{2}}
\section{Applications}
\label{sec:applications}

To illustrate how the group action \bsg{} theorem can be applied, we investigate a few examples.

\subsection{\bsg{} for free actions}

The following theorem generalizes the asymmetric \bsg{} theorem of \cite{tao2010additive} to any group action $G\actson X$ that is \emph{free}, meaning that only the identity element has fixed points.
In particular, it applies to the action of a non-commutative group $G$ on itself by right or left translation.

Recall that a finite subset $S$ of a group is a \emph{$K$-approximate group}\/ if $A=A^{-1}$, $e\in A$, and there is a subset $X\subseteq AA$ of size $|X|=K$ such that $AA\subseteq XA$.
Since $AAA\subseteq XAA\subseteq X^2A$, we have $|A^3|\leq K^2|A|$, so approximate groups have small tripling.
The converse is roughly true.
\begin{lem}[{\cite[Theorem 3.10]{tao2008product}}]
\label{lem:5}
If $A$ is a finite subset of a group such that $|A^3|\leq K|A|$, then $A_{(3)}$ is an $O(K^{15})$-approximate group that contains $A$.
Further, $|A_{(3)}|\ll K^3|A|$.
\end{lem}

\begin{thm}
  \label{thm:4}
Suppose that $G\actson X$ freely.
Let $A\subseteq G$ and $Y\subseteq X$ be finite sets such that $A\subseteq\sym_\alpha(Y)$.

Fix $J>0$ 
and let
\[
L=\alpha^{-2^J} \left( \frac{|Y|}{|A|} \right)^{1/J}.
\]

There exists an $O(L^C)$-approximate group $S$ of size $|S|\ll L^C|Y|$ such that for some element $g$ in $G$
\begin{equation}
  \label{eq:46}
  |A\cap gS| \gg L^{-C}|A|.
\end{equation}
Further, there exist subsets $Y'\subseteq Y$ and $Z\subseteq Y'$ such that
\begin{equation}
  \label{eq:47}
  |Y'|\gg L^{-C}|Y|,
\end{equation}
\begin{equation}
  \label{eq:48}
Y'\subseteq S^{-1}S(Z),
\end{equation}
and
\begin{equation}
  \label{eq:49}
  |Z| \ll L^C\frac{|Y|}{|S|}.
\end{equation}
\end{thm}
\begin{proof}
  By Proposition~\ref{prop:6}, $|\sym_\alpha(Y)|\leq \alpha^{-1}|Y|$, so we may apply Theorem~\ref{thm:3} with
\[
K= \left( \frac{|Y|}{\alpha_J|A|} \right)^{1/J}\approx_J L.
\]
Let $S=(A_*)_{(3)}$; by Lemma~\ref{lem:5}, $S$ is a $O(L^C)$ approximate group containing $A_*$.
Thus by part (2) of Theorem~\ref{thm:3}, there is an element $g$ in $G$ such that
\[
\max \left( |A\cap g^{-1}(gSg^{-1})|, |A\cap gS| \right) \gg L^{-C}|A|.
\]
Further, since $|S|\ll L^C|A_*|$ and $A_*\subseteq g_*\sym_{\alpha_J}(Y)$, we have $|S|\ll L^{C+1}|Y|$.

Since $\fix(g)=\emptyset$ for any non-identity element $g$, by part (3) of Theorem~\ref{thm:3}, there are subsets $Y'\subseteq Y$ and $Z\subseteq Y'$ with $|Y'|\gg L^{-C}|Y|$, $Y'\subseteq S^{-1}S(Y)$, and $|Z|\ll L^C|Y|/|S|$.
\end{proof}

WIth minor modifications, we can prove a more general result.
\begin{thm}[\bsg{} for almost free actions]
  \label{thm:5}
Suppose that $G\actson X$ and that $|\fix(g)|< n$ for all $g\not=e$.

Let $A\subseteq G$ and $Y\subseteq X$ be finite sets such that $A\subseteq\sym_\alpha(Y)$.

Fix $J>0$ and let
\[
L=\alpha^{-n\cdot 2^J} \left( \frac{|Y|}{|A|} \right)^{n/J}.
\]

Then there exists an $O(L^C)$-approximate group $S$ of size $|S|\ll_J \alpha^{-n\cdot 2^J}|Y|^n$, an element $g$ in $G$, subsets $Z\subseteq Y'\subseteq Y$ such that \eqref{eq:46}, \eqref{eq:47}, and \eqref{eq:48} hold, and \eqref{eq:49} is replaced by
\[
|Z| \ll L^C\max \left( \frac{|Y|}{|S|}, \sqrt{n|Y|} \right).
\]
\end{thm}
The proof is the same as that of Theorem~\ref{thm:4}, except that the sum in \eqref{eq:30} is non-empty, so instead we use the bound
\[
\frac 1{|S||Z|}\sum_{g\not=e}|\fix(g)\cap Z|r_{S^{-1}S}(g) < \frac{n|S|}{|Z|}.
\]

\subsection{\bsg{} for linear fractional transformations}
\label{sec:low-dimens-exampl}

In this section, we briefly review some related results from the literature, and indicate how they can be proved using the methods from this paper.
In particular, we would like to emphasize that Theorem~\ref{thm:3} is roughly equivalent to Bourgain and Gamburd's $\ell^2$-flattening strategy \cite{bourgain2008uniform}.

Elekes studied rich affine transformations in \cite{elekes1997linear,elekes1998linear,elekes1999linear,elekes2002versus}, and together with Kiraly, studied rich linear fractional transformations in \cite{elekes2001combinatorics}.
Further work on rich affine transformations was done by Borenstein and Croot \cite{borenstein2010lines} and Amirkhanyan, Bush, Croot, and Pryby \cite{amirkhanyan2017lines}.
The paper \cite{murphy2017upper} extends this later work using Theorem~\ref{thm:3} and shows how to use the tools developed in Section~\ref{sec:group-acti-comb} to recover Elekes' results (as well as extend them to transformations over other fields).

We will show how the tools of this paper can be used to reprove the results of Elekes and Kir\'{a}ly.
For a field $\F$ and $g\in SL_2(\F)$,
\[
g=
\begin{pmatrix}
  a & b\\
  c & d\\
\end{pmatrix},
\]
let $\Gamma_g\subseteq \F\times\F$ denote the curve
\[
cxy - ax+dy-b=0.
\]
\begin{thm}[{\cite[Corollary 35]{elekes2001combinatorics}}]
\label{thm:9}
For $X,Y\subseteq\C$ with $n\leq |X|,|Y|\leq Cn$, and $A\subseteq SL_2(\C)$ with $|A|=n$, if
\[
|\{(x,y,g)\in X\times Y\times A\colon (x,y)\in\Gamma_g\}| \geq \alpha n^2,
\]
then there is an element $g\in SL_2(\C)$ and an abelian subgroup $H\leq SL_2(\C)$ such that
\[
|A\cap gH|\gg (\alpha/C)^{O(1)}n.
\]
\end{thm}
The following weaker analog for finite fields is new.
\begin{thm}
  \label{thm:10}
For $X,Y\subseteq\F_q$ with $n\leq |X|,|Y|\leq Cn$, and $A\subseteq SL_2(\F_q)$ with $|A|=n^3$, if
\[
|\{(x,y,g)\in X\times Y\times A\colon (x,y)\in\Gamma_g\}| \geq \alpha n^4,
\]
then there is an element $g\in SL_2(\F_q)$ and a proper subgroup $H\leq SL_2(\F_q)$ such that
\[
|A\cap gH|\gg (\alpha/C)^{O(1)}n^3.
\]

Further, there is a subset $W$ of $X\cup Y$ of size $|W|\gg (\alpha/C)^{O(1)}n$ that is covered by $O((C/\alpha)^{O(1)}n^{1/2})$ orbits of $H$
\end{thm}

In \cite{bourgain2012modular,shkredov2018asymptotic,moschevitin2018popular} further results on rich linear fractional transformations were obtained using $\ell^2$-flattening.
These results can be recovered by Theorem~\ref{thm:3} as well; we think it is instructive to compare these methods, therefore we will give a short proof of a result for rich linear fractional transformations.
\begin{thm}[{\cite[Proposition 1]{bourgain2012modular}}]
  \label{thm:8}
For all $\epsilon>0$ and $r>1$, there is a $\delta>0$ such that the following holds.
Let $p$ be a large prime.
If $A\subseteq\F_p$ and $S\subseteq SL_2(\F_p)$ satisfy
\begin{enumerate}
\item $1\ll |A|\ll p^{1-\epsilon}$
\item $|A|<|S|^r$
\item $|S\cap gH|<|S|^{1-\epsilon}$ for any proper subgroup $H\leq SL_2(\F_p)$ and any $g\in SL_2(\F_p)$,
\end{enumerate}
then
\[
|\{(x,y,g)\in A\times A\times S\colon (x,y)\in \Gamma_g\}| < 2|A|^{1-\delta}|S|.
\]
\end{thm}
This is analogous to the results on ``rich lines in grids'' proved in \cite{murphy2017upper}.

\begin{lem}
  \label{lem:7}
Suppose $A\subseteq\F_p,S\subseteq SL_2(\F_p)$, and 
\[
|\{(x,y,g)\in A\times A\times S\colon (x,y)\in \Gamma_g\}| \geq 2|A|^{1-\delta}|S|
\]
for some $\delta >0$.
Then there exists a subset $P\subseteq S$ such that $P\subseteq\sym_\alpha(A)$ for $\alpha=|A|^{-\delta}$ and $|P|\geq \alpha|S|$.
\end{lem}
\begin{proof}
As in \cite{murphy2017upper}, we rephrase the problem in terms of symmetry sets.
We have
\[
(x,y)\in \Gamma_g \iff y = gx := \frac{ax+b}{cx+d}.
\]  
Let $P\subseteq S$ denote the set of transformations such that $|A\cap gA|\geq |A|^{1-\delta}$.
By Lemma~\ref{lem:Pop}, we have
\[
\sum_{g\in P}|A\cap gA| \geq \frac 12\sum_{g\in S}|A\cap gA| \geq |A|^{1-\delta}|S|,
\]
hence
\[
|P|\geq |A|^{1-\delta}|S|.
\]
By definition, $P\subseteq\sym_\alpha(A)$, where $\alpha=|A|^{-\delta}$.
\end{proof}

Theorems~\ref{thm:10} and \ref{thm:8} will follow from an auxiliary result.
\begin{lem}
  \label{lem:8}
Given subsets $A\subseteq\F_p$ and $P\subseteq SL_2(\F_p)$, if $P\subseteq\sym_\alpha(A)$, then for any integer $J>0$, we have either
\begin{enumerate}
\item $|A|>(\alpha/2)^{2^J}p$ or $|A| < (\alpha/2)^{-2^J}|Y|^{O(1/J)}$, or
\item there is a proper subgroup $H\leq SL_2(\F_p)$ and an element $g\in SL_2(\F_p)$ such that $|P\cap gH|\gg (\alpha/2)^{O(2^J)}|A|^{-O(1/J)}|P|$.
\end{enumerate}
\end{lem}
The proof of Lemma~\ref{lem:8} requires Helfgott's \emph{product theorem for $SL_2(\F_p)$}\/ \cite{helfgott2008growth}, and a symmetry set bound for $SL_2(\F_p)$ acting on $\P^1(\F_p)$ by linear fractional transformations.
Rudnev and Shkredov \cite{rudnev2018growth} proved a version of the product theorem for $SL_2(\F_p)$ with explicit constants, improving work of Kowalski \cite{kowalski2013explicit}.
\begin{thm}[Growth in $SL_2(\F_p)$]
  \label{thm:9}
  Let $p$ be prime and let $A\subseteq SL_2(\F_p)$ be a set of generators.
  Either
  \begin{equation}
    \label{eq:62}
    (A\cup A^{-1}\cup \{e\})^3=SL_2(\F_p)
  \end{equation}
  or
  \begin{equation}
    \label{eq:63}
    \left( 3\frac{|A^3|}{|A|} \right)^3|A| \geq|(A\cup A^{-1}\cup \{e\})^3|\geq |A|^{1+\delta}
  \end{equation}
  where $\delta=\frac 1{3024}$.
\end{thm}
Note also that
\begin{equation}
  \label{eq:67}
      \left( 3\frac{|A^3|}{|A|} \right)^3|A| \geq|(A\cup A^{-1}\cup \{e\})^3|.
\end{equation}


\begin{proof}[Proof of Lemma~\ref{lem:8}]
    Since $PSL_2(\F)$ embeds into $PGL_2(\F)$, which acts simply 3-transitively on $\P^1(\F)$, no element of $PSL_2(\F)$ has more than two fixed points, except the identity.
By Theorem~\ref{thm:5}, there is an $O(L^C)$-approximate group $S$, where $L=\alpha^{-2^J} (|Y|/|A|)^{1/J}$, such that
\[
|A\cap gS|\gg L^{-C}|A|
\]
for some $g\in PSL_2(\F)$.
Further, there exists $Y'\subseteq Y$ and $Z\subseteq Y'$ such that $|Y'|\gg L^{-C}|Y|$, $Y'\subseteq S^{-1}S(Z)$, and $|Z|\ll L^C\max(|Y|/|S|, \sqrt{|Y|})$.

By Theorem~\ref{thm:9}, either $S^3=PSL_2(\F)$, $|S^3|\geq |S|^{1+\delta}$, or $S$ is contained in a proper subgroup.
Since $|S^3|\ll L^{C}|S|$ and $|S|\ll_J \alpha^{-3\cdot 2^J}|Y|^3$, if $|S|\gg L^{C/\delta}$ and $|Y| \ll_J \alpha^{2^J}L^Cp$, then $S$ is contained in a proper subgroup $H$.

If $|S|\ll L^{C/\delta}$, then $|A|\ll L^{C'}$, so $|A|\ll \alpha^{-C'\cdot 2^J}|Y|^{C'/J}$
\end{proof}
Lemma~\ref{lem:8} can also be proved using $\ell^2$ flattening \cite{moschevitin2018popular}.

  



\subsection{Linear actions over fields of characteristic zero}
\label{sec:linear-actions-C}

The following theorem generalizes the work of Elekes \cite{elekes1997linear} and Elekes and Kiraly \cite{elekes2001combinatorics}.
\begin{thm}
  \label{thm:rich-linear-char-zero}
Let $k$ be a field of characteristic zero, let $A$ be a finite subset of $GL_n(k)$, and let $Y$ be a finite subset of $k^n$.

Fix an integer $J>0$ and $1/|Y|<\rho\leq\alpha\leq 1$ such that $\rho \leq (\alpha/2)^{2^J}$, and suppose that for all proper linear subspaces $W\leq k^n$, we have $|Y\cap W|\leq\rho|Y|$.
If $A\subseteq\sym_\alpha(Y)$, then for all integers $J>0$ there is an element $g\in GL_n(k)$ and a nilpotent subgroup $N\leq GL_n(k)$ of step at most $n-1$ such that
\[
|A\cap gN| \gg (\alpha/2)^{O_n( 2^J)}|Y|^{-O_n(1/J)}|A|.
\]
\end{thm}
The assumption that $Y$ does not concentrate in proper subspaces is necessary without further assumptions on $A$, since if $|Y\cap W|\geq\alpha|Y|$ for some subspace $W\leq V$ of codimension at least 2, then the set of transformations stabilizing $W$ is contained in $\sym_\alpha(Y)$, while $\stab(W)$ is not solvable, hence not nilpotent.
However, the relationship between $\rho$ and $\alpha$ in Theorem~\ref{thm:rich-linear-char-zero} may not be optimal.

The proof of Theorem~\ref{thm:rich-linear-char-zero} requires the following product theorem in $GL_n(k)$, where $k$ is a field of characteristic zero, is due to Breuillard, Green, and Tao \cite[Theorem 2.5]{breuillard2011approximate}.
\begin{thm}
  \label{thm:BGT-char-zero}
Let $k$ be a field of characteristic zero.
Suppose that $A\subseteq GL_n(k)$ is a finite subset such that $|A^3|\leq K|A|$.
Then there is a subset $B\subseteq GL_n(k)$, an element $g\in GL_n(k)$, and a constant $C$ depending on $n$ such that
\begin{itemize}
\item $|B|\ll K^C|A|$,
\item $|A\cap gB|\gg K^{-C}|A|$,
\item $|B^3|\ll K^C|B|$, and
\item $B$ generates a nilpotent group of step at most $n-1$.
\end{itemize}
\end{thm}

We also require the following symmetry set bound, which is used in the proof of Theorem~\ref{thm:rich-linear-char-zero}, whose proof we defer to the end of the section.
\begin{prop}[Bound for linear actions]
  \label{prop:11}
Suppose that $G$ acts linearly on an $n$-dimensional vector space $V$.
Let $Y\subseteq V$ be a finite subset and suppose that there are parameters $1/|Y|<\rho<\alpha\leq 1$ such that for any subspace $W\leq V$, we have $|Y\cap W|\leq \rho|Y|$.
If the action of $G$ on $V$ is \emph{faithful}, then
\[
|\sym_\alpha(Y)|\leq\frac{|Y|^n}{(\alpha-1/|Y|)(\alpha-\rho)^{n-1}}.
\]
\end{prop}
For instance, if $\rho\leq\alpha/2$, we have
\[
|\sym_\alpha(Y)|\ll \frac{|Y|^n}{\alpha^n}.
\]
The requirement for $G$ to act faithfully on $V$ simply means that $G$ can be realized as a subgroup of $GL(V)$.
\begin{remark}
Stronger bounds than Proposition~\ref{prop:11}, with stricter hypotheses, have appeared before.
Solymosi and Tao~\cite{solymosi2011incidence} proved that if $G$ acts on a hyperplane $H$ of $V$ by affine transformations, where $V$ is a real finite dimensional vector space, then for $A\subseteq G$ and a finite set $Y\subseteq H$ we have
\begin{equation}
  \label{eq:61}
  |A\cap\sym_\alpha(Y)|\ll_\epsilon \frac{|Y|^{1+\epsilon}}{\alpha^3},
\end{equation}
provided that no two elements of $A$ have a common fixed point.
Elekes conjectured \cite{elekes2011dimension} that \eqref{eq:61} should hold (perhaps with higher powers of $\alpha$) with $\epsilon=0$ for a variety of algebraic actions under certain hypotheses.
In Theorem 2.3 of the same paper, Elekes proved such a result (without explicit dependence on $\alpha$) for affine transformations of $\R^2$ acting on sets $Y$ that are \emph{proper 2-dimensional}, meaning that they can be cut into singletons by $O(|Y|^{1/2})$ lines.
Guth and Katz~\cite{guth2015erdos} proved a bound of the form \eqref{eq:61} with $\epsilon=0$ for isometries of $\R^2$, confirming Conjecture 2.1 of \cite{elekes2011dimension}, and as a consequence, settling the Erd\H{o}s distinct distance problem, up to logarithmic factors.

Though Proposition~\ref{prop:11} is weaker than the bounds cited above, it is sufficiently strong to prove Theorem~\ref{thm:rich-linear-char-zero}.
This illustrates that the group action \bsg{} theorem only requires bounds on $|\sym_\alpha(Y)|$ that are \emph{polynomial}\/ in $|Y|$ and $\alpha$, whereas the previous methods of Elekes (and Elekes and Kiraly) require bounds that are \emph{linear}\/ in $|Y|$.
\end{remark}



The proof of Proposition~\ref{prop:11} is based on the following variation of Proposition~\ref{prop:6}.
\begin{lem}
  \label{lem:3-upper}
Suppose that $G\actson X$.
For a positive integer $n$, we have the diagonal action $G\actson X^n$ given by $g(x_1,\ldots,x_n)=(gx_1,\ldots,gx_n)$.

Given $Y\subseteq X$, suppose there is a subset $Y_*\subseteq Y^n$ such that for all $g\in\sym^{G\actson X}_\alpha(Y)$, we have $|Y_*\cap gY^n|\geq \alpha_*|Y|^n$ and for any $(y_1,\ldots,y_n)\in Y_*$ we have $|\stab(y_1,\ldots,y_n)|\leq C$.
Then
\[
|\sym_\alpha(Y)|\leq \frac{C|Y_*|}{\alpha_*}\leq \frac{C|Y|^n}{\alpha_*}.
\]
\end{lem}
\begin{proof}
By the Elekes-Sharir paradigm (\ref{eq:elekes-sharir}), we have
\begin{align*}
  \alpha_*|Y|^n|\sym^{G\actson X}_\alpha(Y)| &\leq \sum_{g\in\sym_\alpha(Y)}|gY_*\cap Y^n|\\
& = \sum_{\vec y_1\in Y_*,\vec y_2\in Y^n}|\sym_\alpha(Y)\cap \trans(\vec y_1,\vec y_2)|.
\end{align*}
On the other hand, since $\vec y_1\in Y_*$ we have
\[
|\sym_\alpha(Y)\cap \trans(\vec y_1,\vec y_2)|\leq |\trans(\vec y_1,\vec y_2)| = |\stab(\vec y_1)|\leq C.
\]
\end{proof}

If $G\actson X$ and $x_1,\ldots,x_n\in X$, then we will use $\stab(x_1,\ldots,x_n)$ to denote the \emph{point-wise stabilizer}\/ of the set $\{x_1,\ldots,x_n\}$.
That is,
\[
\stab(x_1,\ldots,x_n)=\bigcap_{i=1}^n \stab(x_i).
\]
If $\vec x = (x_1,\ldots,x_n)$, we write $\stab(\vec x)$ for $\stab(x_1,\ldots,x_n)$.

\begin{proof}[Proof of Proposition~\ref{prop:11}]
We can consider $Y^n$ as a set of $n\times n$ matrices; let $Y_*=Y^n\cap GL(V)$.
Then for any $\vec y\in Y_*$, we have $\stab(\vec y)=\{e\}$, since the action of $G$ on $V$ is \emph{faithful}.

Suppose that $g\in\sym^{G\actson V}_\alpha(Y)$.
We wish to show that $|gY_*\cap Y^n|\geq \alpha_*|Y|^n$ for $\alpha_*=(\alpha-1/|Y|)(\alpha-\rho)^{n-1}$.
That is, we wish to count the number of tuples $(y_1,\ldots,y_n)\in Y^n$ such that $gy_i\in Y$ for $i=1,\ldots,n$ and such that $y_{i+1}$ is not in the subspace spanned by $y_1,\ldots, y_i$.
There are at least $\alpha|Y|-1$ non-zero choices for $y_1$.
Given $y_1$, there are at least $(\alpha-\rho)|Y|$ choices of $y_2$ such that $gy_2\in Y$ and $y_2$ is not in the line determined by $y_1$.
Similarly, there are at least $(\alpha-\rho)|Y|$ choices of $y_3$ such that $gy_3\in Y$ and $y_3$ is not in the plane spanned by $y_1$ and $y_2$, and so on.

The desired bound on $|\sym_\alpha(Y)|$ now follows from Lemma~\ref{lem:3-upper} with $C=1$ and $\alpha_*=(\alpha-1/|Y|)(\alpha-\rho)^{n-1}$.
\end{proof}

Now we are ready to prove the theorem.
\begin{proof}[Proof of Theorem~\ref{thm:rich-linear-char-zero}]
  Since $|W\cap Y|\leq\rho|Y|$ for all proper subspaces $W\leq V$, by Proposition~\ref{prop:11} we have
\[
|\sym_\alpha(Y)|\leq \frac{|Y|^n}{(\alpha-1/|Y|)(\alpha-\rho)^{n-1}}
\]
for all $1/|Y|\leq\rho<\alpha\leq 1$.

Applying Theorem~\ref{thm:3}, we find an element $g_*\in GL_n(k)$ and a subset $A_*\subseteq g_*\sym_{\alpha_J}(Y)$ with
\[
|A_*^3|\ll \left( \frac K{\alpha_J} \right)^C |A_*|,
\]
where
\[
  K= \left( \frac{|\sym_{\alpha_J}(Y)|}{|A|} \right)^{1/J} \ll \left( \frac{|Y|^n}{\alpha_J(\alpha_J-\rho)^{n-1}} \right)^{1/J} \ll \left( \frac{|Y|^n}{\alpha_J^n} \right)^{1/J},
\]
since $\rho\leq (\alpha/2)^{2^T}=\frac 12\alpha_J$.

By Theorem~\ref{thm:BGT-char-zero}, there is a subset $B\subseteq GL_n(k)$ and an element $g\in GL_n(k)$ such that $B$ generates a nilpotent group $N$ of step at most $n-1$ and
\[
|A_*\cap gB|\gg \left( \frac{\alpha_J}K \right)^{C_n}|A_*|,
\]
where $C_n>0$ is a constant depending only on $n$.

By part (2) of Theorem~\ref{thm:3}, there is an element $g'$ such that
\[
|A^{-1}\cap g'gB|+|A\cap g'gB| \gg \left( \frac{\alpha_J}K \right)^{O(C_n)}|A|.
\]
It follows that
\[
|A\cap g_0N|+|A\cap g_0N'|\gg \left( \frac{\alpha_J}{|Y|^{1/J}} \right)^{O_n(1)}|A|,
\]
where $N'$ is a conjugate of $N$, hence nilpotent.
Noting that $\alpha_J=2(\alpha/2)^{2^J}$ completes the proof.
\end{proof}




\appendix

\section{Covering lemmas}
\label{sec:covering-lemmas}

To prove our covering lemmas, we need the orbit-stabilizer theorem \emph{for sets}\/ (rather than for groups) \cite{helfgott2015growth}.
\begin{lem}
\label{lem:4}
Suppose $G\actson X$, $x\in X$, and $A\subseteq G$ is finite.
Then there exists $a_0$ in $A$ such that
\begin{equation}
  \label{eq:x1}
  |(a_0^{-1}A)\cap\stab(x)|\geq\frac{|A|}{|A(x)|},
\end{equation}
and for all finite sets $B\subseteq G$,
\begin{equation}
  \label{eq:x2}
  |BA|\geq |A\cap\stab(x)||B(x)|.
\end{equation}
\end{lem}
Lemma~\ref{lem:4} follows from Corollaries~\ref{cor:growth-in-a-subgroup} and \ref{cor:1}.

Recall the statement of Proposition~\ref{prop:2}.
\begin{repprop}{prop:2}[Ruzsa covering lemma for group actions]
Suppose that \(|A(Y)|=K|Y|\)
Then there exists \(Z\subseteq Y\) such that
\begin{enumerate}
\item \(Y\subseteq A^{-1}A(Z)\)
\item \(|A(Z)|=\sum_{z\in Z}|A(z)|\)
\item \(\displaystyle |Z|\leq  \frac{K|Y|}{|A|}\left(1+\frac 1{|A||Z|}\sum_{g\not=e} |\fix(g)\cap Z|r_{A^{-1}A}(g)\right)\).
\end{enumerate}
\end{repprop}
By a pigeonholing argument, the upper bound 3. implies
\[
|Z|\leq \frac{K|Y|}{|A|}\max_{z\in Z}|A^{-1}A\cap\stab(z)|.
\]
We will give a direct proof of this bound before proving the more complicated upper bound.

The following product construction shows that the upper bound in part 3 is sharp in general.
\begin{ex}
\label{ex:covering-sharp}
Let $G= G_1\times G_2$ act on $\{e\}\times G_2$ by left translation, let $A=G_1\times A_2$, and let $Y=\{e\}\times Y_2$ where $|A_2(Y_2)|\leq K|Y_2|$.
Since $A_2(Y_2)$ is just a product set, by standard Ruzsa covering lemma (or Proposition~\ref{prop:2}, noting that point stabilizers are trivial for $G_2\actson G_2$), there is a subset $Z_2\subseteq Y_2$ such that $Y_2\subseteq A_2^{-1}A_2(Z_2)$ and $|Z_2|\leq K|Y_2|/|Z_2|$.
Setting $Z=\{e\}\times Z_2$, we have $Y\subseteq A^{-1}A(Z)$, and 
\[
|Z|=|Z_2|\leq \frac{K|Y_2|}{|Z_2|} = \frac{K|Y|}{|Z|}|G_1|=\frac{K|Y|}{|Z|}\max_{z\in Z}|A^{-1}A\cap \stab(z)|,
\]
where the last equality follows from $G_1\times \{e\}\subseteq A^{-1}A$.
\end{ex}

\begin{proof}[Proof of Proposition~\ref{prop:2}]
Let \(Z\subseteq Y\) be a maximal subset such that \ref{disjoint} holds; that is, such that the images (``approximate orbits'') \(A(z)\) with \(z\) in \(Z\) are disjoint.

If \(y\in Y\), then by maximality, there exists \(z\) in \(Z\) such that \(A(y)\cap A(z)\) is non-empty.
Thus \(y\in A^{-1}A(Z)\), which proves \ref{cover}.

To prove the last statement, note that by definition of $Z$
\[
\sum_{z\in Z}|A(z)| = |A(Z)| \leq |A(Y)|\leq K|Y|.
\]
On the other hand, by Lemma~\ref{lem:4}, for all $z$ in $Z$ we have
\[
|A(z)|\geq \frac{|A|}{|A^{-1}A\cap\stab(z)|}.
\]
Thus
\[
\sum_{z\in Z}\frac 1{|A^{-1}A\cap\stab(z)|} \leq \frac{K|Y|}{|A|},
\]
which yields
\[
|Z|\leq \frac{K|Y|}{|A|}\max_{z\in Z}|A^{-1}A\cap\stab(z)|.
\]

To prove the upper bound stated in 3. we will use action energy.
By Cauchy-Schwarz we have
\[
|A|^2|Z|^2\leq |A(Z)| E(A,Z) \leq K|Y| E(A,Z).
\]
Since $A(z)\cap A(z')=\emptyset$ for $z\not=z'$, we have
\[
E(A,Z)=\sum_{z\in Z}E(A,\{z\}).
\]
By Proposition~\ref{prop:8},
\[
E(A,\{z\})\leq \sum_{g\in\stab(z)}r_{A^{-1}A}(g).
\]
Hence
\[
|Z|^2\leq \frac{K|Y|}{|A|^2}\sum_{z\in Z}\sum_{g\in\stab(z)}r_{A^{-1}A}(g).
\]
As mentioned in Section~\ref{sec:image-sets}, 
\[
\sum_{z\in Z}\sum_{g\in\stab(z)}r_{A^{-1}A}(g)= |A||Z| + \sum_{g\not=e}|\fix(g)\cap Z|r_{A^{-1}A}(g),
\]
which completes the proof.
\end{proof}

The following proposition is a covering result, like Proposition~\ref{prop:2}, but with weaker hypotheses and conclusion.
Namely, we only assume that $B$ is a set of $\alpha$-rich transformations of $Y$, rather than assuming that the image set $B(Y)$ is small.
As a consequence, we are only able to cover a subset $Y'\subseteq Y$; if $BB^{-1}$ is small, then $Y'$ is large.
\begin{repprop}{prop:8}
Suppose $B\subseteq\sym_\alpha(Y)$.
The there exist subsets $Y'\subseteq Y$ and $Z\subseteq Z'$ such that
\begin{equation*}
\tag{\ref{eq:1}}
  Y'\subseteq B^{-1}B(Z),
\end{equation*}
\begin{equation*}
\tag{\ref{eq:2}}
  |Y'| \geq \frac{\alpha|B|}{2|BB^{-1}|}|Y|,
\end{equation*}
and
\begin{equation*}
\tag{\ref{eq:3}}
|Z|\leq \frac{2|Y|}{\alpha|B|}\max_{z\in Z}|B^{-1}B\cap\stab(z)|.
\end{equation*}
In addition,
\begin{align*}
  |Z|&\leq \frac{4|Y|}{\alpha^2|B|^2}\frac 1{|Z|}\sum_{z\in Z}\sum_{g\in\stab(z)}r_{B^{-1}B}(g) \\
&= \frac{4|Y|}{\alpha^2|B|} \left( 1+\frac 1{|B||Z|}\sum_{g\not=e}|\fix(g)\cap Z|\, r_{B^{-1}B}(g) \right).
\end{align*}  
\end{repprop}

\begin{proof}[Proof of Proposition~\ref{prop:5}]
Since $B\subseteq\sym_\alpha(Y)$, we have
\[
\sum_{g\in B}|Y_g|\geq \alpha|B||Y|.
\]
On the other hand,
\begin{align*}
  \sum_{g\in B}|Y_g|
&= |\{(g,y,y')\in B\times Y\times Y\colon gy=y'\}|\\
&=  \sum_{y\in Y} \left( \sum_{y'\in Y} |B\cap\trans(y,y')| \right).
\end{align*}
Since $\trans(y,y')=g\stab(y)$ for any $g$ in $\trans(y,y')$, if $B\cap\trans(y,y')$ is non-empty, then there is an element $b_{y'}$ in $B$ such that
\[
\trans(y,y')=b_{y'}\stab(y).
\]
Hence $|B\cap\trans(y,y')|=|B\cap b_{y'}\stab(y)|$, which yields
\[
 \sum_{y\in Y} \left( \sum_{y'\in Y} |B\cap b_{y'}\stab(y)|\right) = \sum_{g\in B}|Y_g|\geq\alpha|B||Y|.
\]

Now by Lemma~\ref{lem:Pop} there is a subset $P\subseteq Y$ such that
\[
\sum_{y'\in Y} |B\cap b_{y'}\stab(y)|\geq\frac{\alpha}2|B|
\]
for all $y$ in $P$ and
\begin{equation}
  \label{eq:54}
   \sum_{y\in P} \left( \sum_{y'\in Y} |B\cap b_{y'}\stab(y)|\right) = \sum_{g\in B}|Y_g|\geq\frac{\alpha}2 |B||Y|.
\end{equation}

By Lemma~\ref{lem:4} we have
\[
|B\cap b_{y'}\stab(y)|=|B^{-1}b_{y'}\cap\stab(y)|\leq \frac{|BB^{-1}b_{y'}|}{|B(y)|}\leq \frac{|BB^{-1}|}{|B(y)|}.
\]
Since $B\cap\trans(y,y')$ is non-empty whenever $y'\in B(y)$, we have
\begin{equation}
  \label{eq:55}
\sum_{y'\in Y} |B\cap b_{y'}\stab(y)| \leq |BB^{-1}|\frac{|B(y)\cap Y|}{|B(y)|}\leq |BB^{-1}|.
\end{equation}
By equations~(\ref{eq:54}) and (\ref{eq:55})
\[
|P| \geq \frac{\alpha|B|}{2|BB^{-1}|}|Y|.
\]

The rest of the proof follows the last step of the proof of \cite[Theorem 2.35]{tao2010additive} closely, imitating the proof of Ruzsa's covering lemma (see the proof of Proposition~\ref{prop:2}).
Choose $Z\subseteq P$ to be a maximal subset such that
\begin{equation}
  \label{eq:56}
  |B(Z)\cap Y|=\sum_{z\in Z}|B(z)\cap Y|.
\end{equation}
Let $Y'=P$.
Then if $y\in Y'$, we have $B(y)\cap B(z)\not=\emptyset$ for some $z$ in $Z$, so
\[
Y'\subseteq B^{-1}B(Z).
\]

To give an upper bound for $|Z|$, we will double count the size of the set $E:=\{(b,z)\in B\times Z\colon b(z)\in Y\}$.
Since
\[
|E|=\sum_{y\in Y}r_{B(Z)}(y) = \sum_{y\in Y}\sum_{z\in Z}|B\cap\trans(z,y)|
\]
and $Z\subseteq P$, by definition of $P$ we have
\[
|E| =\sum_{z\in Z} \left( \sum_{y\in Y}|B\cap\trans(z,y)| \right)\geq \frac{\alpha|B||Z|}2.
\]
On the other hand,
\[
|E| = \sum_{y\in Y} \left( \sum_{z\in Z}|B\cap\trans(z,y)| \right)\leq |B(Z)\cap Y| \sup_{y\in Y}\sum_{z\in Z}|B\cap\trans(z,y)|.
\]
Since $B(z)\cap B(z')=\empty$ for all distinct pairs $z,z'\in Z$, for all $y\in Y$ such that the inner sum on the right-hand side of the above equation is non-zero, there is a $z_y\in Z$ such that
\[
\sum_{z\in Z}|B\cap\trans(z,y)|=|B\cap\trans(z_y,y)|.
\]
Since $\trans(z_y,y)=g_y\stab(z)$ for some $g_y\in B$, we have
\[
|E| \leq |B(Z)\cap Y|\sup_{z\in Z, b\in B}|b^{-1}B\cap\stab(z)|\leq |Y|\sup_{z\in Z}|B^{-1}B\cap\stab(z)|.
\]
All together, we have
\[
|Z| \leq \frac{2|Y|}{\alpha|B|}\,\sup_{z\in Z}|B^{-1}B\cap\stab(z)|.
\]

We can give an alternate bound with the supremum replaced by an average.
Let $E=\{(b,z)\in B\times Z\colon b(z)\in Y\}$, as above; we know that $|E|\geq \alpha|B||Z|/2$.
Then $B_E(Z)=B(Z)\cap Y$ and
\[
|E|=\sum_{y\in B(Z)\cap Y}r_{B(Z)}(y).
\]
By Cauchy-Schwarz, we have
\[
|E|^2|\leq |B(Z)\cap Y| \sum_{y\in Y}r_{B(Z)}^2(y).
\]
Since $bz=b'z'=y$ implies that $z=z'$, we have
\begin{align*}
  \sum_{y\in Y}r_{B(Z)}^2(y) &= |\{(b,b',z,z',y)\in B\times B\times Z\times Z\times Y\colon bz=bz'=y\}|\\
&\leq \sum_{z\in Z}\sum_{g\in\stab(z)}r_{B^{-1}B}(g).
\end{align*}
The last inequality follows from a similar argument to the proof of Corollary~\ref{cor:1}.

Combining these equations yields
\[
\frac{\alpha^2|B|^2|Z|^2}4\leq |B(Z)\cap Y|\sum_{z\in Z}\sum_{g\in\stab(z)}r_{B^{-1}B}(g).
\]
Thus
\[
|Z|^2 \leq \frac{4|Y|}{\alpha^2|B|^2}\sum_{z\in Z}\sum_{g\in\stab(z)}r_{B^{-1}B}(g).
\]
This can also be rearranged as
\[
|Z|^2 \leq \frac{4|Y|}{\alpha^2|B|^2} \left( |B||Z| + \sum_{g\not=e}|\fix(g)\cap Z|r_{B^{-1}B}(g) \right),
\]
or
\[
|Z| \leq \frac{4|Y|}{\alpha^2|B|} \left(1+ \frac 1{|B||Z|} \sum_{g\not=e}|\fix(g)\cap Z|r_{B^{-1}B}(g) \right).
\]
\end{proof}

\bibliographystyle{plain}
\bibliography{/Users/brendan/Dropbox/library}

\end{document}